\documentclass[preprint,10pt,authoryear]{elsarticle}

\newcommand{\cK}{\mathcal{K}}

\newcommand{\CA}[1]{\textcolor{red}{#1}}

\usepackage{amsmath} 
\usepackage{amssymb}
\usepackage{pifont}
\usepackage[table,xcdraw]{xcolor}
\usepackage{rotating}
\usepackage{multirow}
\usepackage{comment}
\usepackage{amssymb}
\usepackage[inline,shortlabels]{enumitem}

\usepackage{xcolor}
\usepackage[hyphens]{url}
\usepackage{hyperref}

\usepackage[normalem]{ulem}
\hypersetup{colorlinks=true,breaklinks=true} 
\usepackage{setspace}
\usepackage{cancel}
\begin{document}

\begin{frontmatter}



\title{Emerging Optimization Problems for Distribution in Same-day Delivery}


\author{Yuanyuan Li\corref{cor1}}
\cortext[cor1]{Corresponding author}
\ead{yuanyuan.li@essec.edu} 

\author{Claudia Archetti}
\ead{archetti@essec.edu} 

\author{Ivana Ljubi\'c}\ead{ivana.ljubic@essec.edu}

\address{IDS Department, ESSEC Business School, Cergy-Pontoise, France}

\begin{abstract}
Same-day deliveries (SDD) have become a new standard to satisfy the ``instant gratification'' of online customers. 
Despite the existing powerful technologies deployed in last-mile delivery, 
SDD services face new 
decision-making challenges related to the 
trade-off between delivery cost and time. In addition, new challenges related to environmental issues, customer satisfaction, or fairness arise. 
Researchers have explored various approaches to face these challenges in the context of SDD, where stochastic and dynamic data uncertainty plays a fundamental role.
In this paper, we carefully review the emerging routing problems and solutions proposed in the existing literature for SDD services. We address the questions related to how to deal with dynamic arrival times of orders, how to allocate time slots to deliveries, how to select the right delivery options, how to design pickup and delivery routes, or how to partition the delivery areas and decide the composition of the fleet. We also formulate and compare models for representative problems elaborating on the pros and cons that might guide practitioners in choosing the most appropriate objectives and constraints. Finally, we sketch challenges and identify future research directions.
\end{abstract}


\end{frontmatter}


\allowdisplaybreaks 
\section{Introduction}
\label{intro}

The COVID-19 pandemic was responsible for the most drastic change in consumption patterns in recent years. House-bound consumers shifted their purchases to online shopping, sometimes exceeding the pre-pandemic levels in what can be termed a ``revenge purchase'' \citep{park2022changes}. 
Powerful technologies have emerged: autonomous vehicles and drones are now increasingly used for shipping, and IoT-enabled devices are used for tracking deliveries. 
This combination of changes in consumer behavior and emerging of new technologies 
has not only affected
the volume of online shopping but has also raised the consumers' expectations with respect to shorter delivery times\footnote{https://www.statista.com/statistics/1274950/global-online-shoppers-wished-changes-on-delivery/}.
Conventional delivery (the delivery of parcels usually in 3 to 5 days from their order time) is increasingly being replaced by \emph{same-day delivery} (SDD, the delivery of goods not longer than one day from their purchase) or even more extreme \emph{instant delivery} (an on-demand delivery within two hours from order request) \citep{dablanc2017rise}). It is estimated that by 2025 same-day and instant deliveries will reach 20 to 25 percent of the market share. They can grow substantially further, notably, if the service is expanded to rural areas \citep{joerss2016customer}. 

Although SDD services bring lots of conveniences to consumers, they also create a series of challenges for providers, as the arrival of requests and their acceptance, picking and packing, and delivery to customer locations all need to happen within the same day. This involves restructuring the existing delivery infrastructures, ranging from warehouse design to delivery processes. 
Recent surveys 
by \cite{bock2021pro}, \citet{archetti2021recent} and \citet{boysen2021last} {briefly review some strategies in online routing designed for urgent delivery processes and} cover a broad range of relevant and emerging optimization problems in e-commerce and last-mile delivery. However, they do not focus on SDD specifically and on the corresponding optimization problems. 
The current work aims at filling this gap by providing a critical review of the delivery challenges and corresponding solutions in SDD services. 
Due to the limited space, we mainly focus on distribution operations and do not cover problems related to, e.g., inventory management. 

\begin{figure}[!h]  \centering\includegraphics[width=0.9\textwidth]{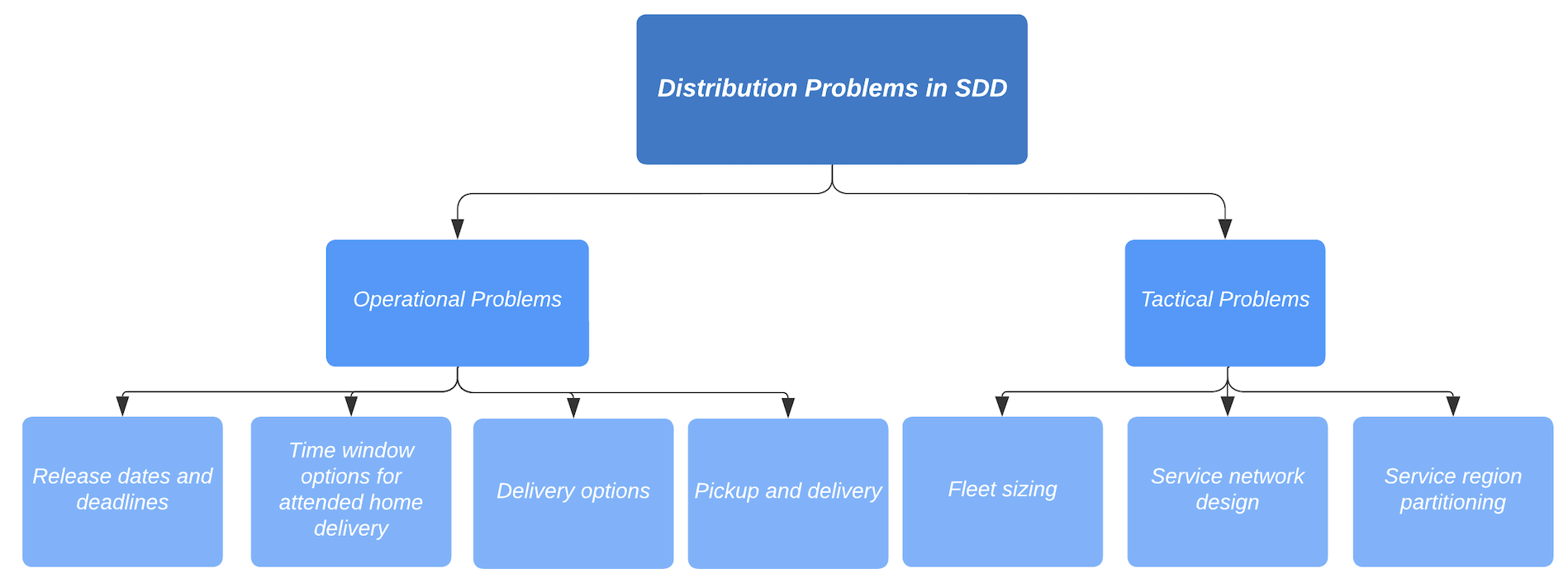}
    \caption{Classification scheme of SDD problems}
    \label{fig:SDDprob}
\end{figure}

The main aim of this paper is to analyze the emerging optimization problems in SDD services. Specifically, we focus on problems that embed new features related to the nature of SDD, and which, as such, render the corresponding optimization problems different from the ones faced in the former literature. As mentioned above, we focus on delivery operations and, more specifically, the routing problems related to the distribution process. In SDD, order arrival, order assembly, vehicle routing and delivery all happen within a single day. Since orders arrive during distributions, we focus on \textit{order-based fulfillment systems} as defined by \cite{WAMUTH2023801}. Notice that these systems work differently than the \textit{periodic fulfillment} ones, where a cut-off time is set to order arrivals and distribution happens afterwards, and more flexibility is provided in dealing with dynamic order arrivals. This paper focuses on the operational and tactical planning levels of SDD services (as shown in Figure \ref{fig:SDDprob}). We acknowledge that the lists of operational and tactical problems presented herein may not be exhaustive, as the landscape of SDD problems is continually evolving. The following classes of relevant routing problems in SDD have been identified: 
\begin{itemize}
    \item Routing problems with release dates and deadlines (see Section \ref{RD_deadlines});
    \item Time window options in attended home delivery services (see Section \ref{time_window});
    \item Delivery options including where to deliver and how to deliver (see Section \ref{delivery});
    \item Pickup and delivery problems (see Section \ref{pickup_delivery}); 
    \item Districting problems involving fleet size and mix vehicle routing problem (see Section \ref{Districting}).
\end{itemize}   

In the five classes of problems we present in this work, the first four mainly deal with operational decisions, following the process we sketch in Figure \ref{fig:SDDprocess}; the last one (districting) mainly concerns tactical decisions such as fleet sizing, service network design, and service region partitioning. In operational problems, the ones with release dates and deadlines serve as the basis for many other SDD problems that often include additional complexities beyond these temporal constraints. Thus we start our discussion in the subsequent sections by introducing the problems with release dates and deadlines.

We present mathematical formulations for routing problems with release dates, time window options, and delivery options with pickup stations. The models are built for basic settings, i.e., deterministic case and a single vehicle. As explained in the following, these models help the understanding of each problem setting and can be used as a basis of solution approaches for more complex systems or to design \textit{a priori} solutions (\cite{klapp2018dynamic}), for competitive analysis in the case of stochastic settings (\cite{li2022reinforcement}, \cite{stroh2022tactical}), or as tools embedded in solution approaches for more general problems. 
For the other two classes of problems, i.e., pickup \& delivery and districting, we do not present problem formulations as they can be obtained by combining the traditional models with the ones introduced in Formulation \eqref{eq:base1}, i.e., the first formulation presented in the current work. 
In addition, we compare the formulations proposed on example instances with the aim of illustrating the pros and cons. 
%
Finally, we outline challenges and specify future research directions. 

The paper is structured as follows. Section \ref{overview} provides an overview of the SDD problems by structuring process flows, identifying the main sources of uncertainty and approaches to handle the uncertainties. Subsequently, it presents classification schemes for distribution problems in SDD. The different classes of operational and tactical problems are reviewed in Sections \ref{RD_deadlines}-\ref{Districting}. Section \ref{model-comparison} compares the formulations and the corresponding solutions on small instances. Section \ref{challenges} presents challenges and future research directions. Section \ref{conclusion} concludes the paper.



\section{Overview on the SDD Problems}
\label{overview}

In this section, we first highlight the primary stakeholders in SDD services, we explain the main process flows, and discuss the various uncertainty factors in Section \ref{process_uncertainty}. 
Then we present a classification scheme for SDD problems in Section \ref{classification_scheme}.

\subsection{Major stakeholders in SDD, process flows and sources of uncertainty}
\label{process_uncertainty}

\begin{figure}[!h]  \centering\includegraphics[width=0.6\textwidth]{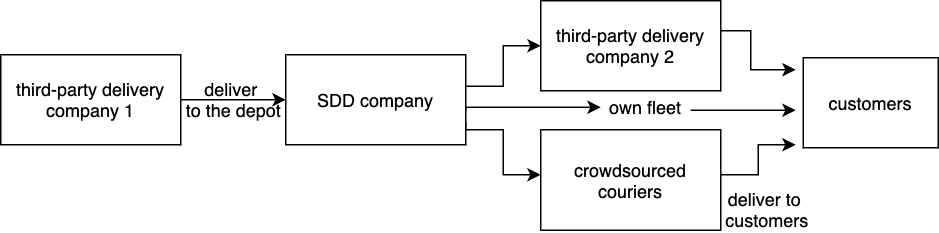}
    \caption{SDD stakeholders}
    \label{fig:stakeholders}
\end{figure}
\begin{figure}[!h]  \centering\includegraphics[width=0.65\textwidth]{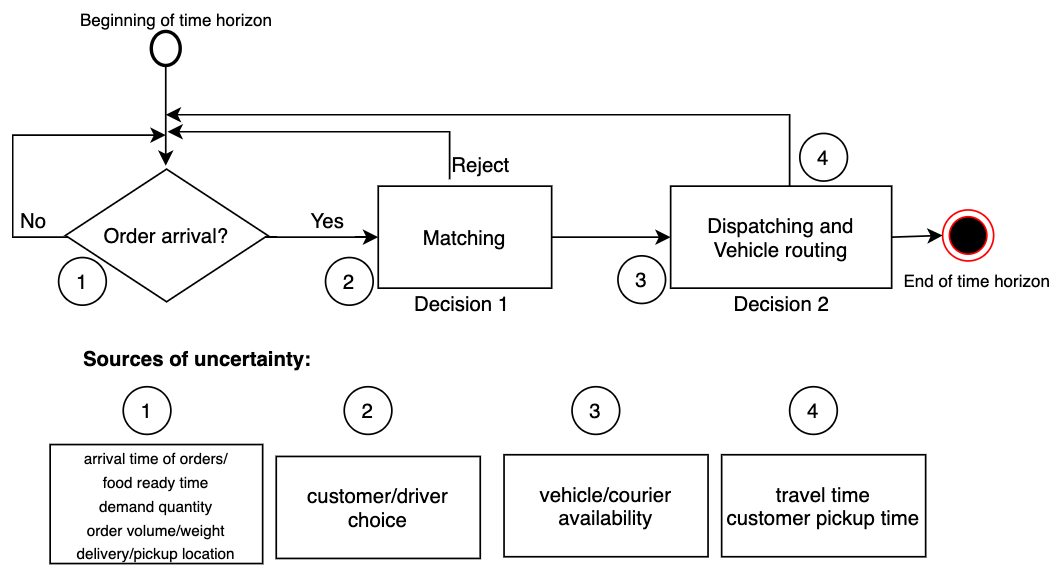}
    \caption{SDD process with sources of uncertainty}
    \label{fig:SDDprocess}
\end{figure}
SDD services manage day-to-day operations in an environment characterized by continuous customer order arrivals. The products being delivered through this system can be newspapers\footnote{https://www.newspaperdirect-asia.com/}, groceries, medicines, and fashion-related items, all of which are frequently requested through SDD services.

The principal stakeholders involved in this distribution process include retail enterprises engaged in online product sales and provision of SDD, customers who request the service, third-party delivery entities responsible for delivering parcels to depots, and/or from depots to customers, and crowdsourced couriers who deliver using their personal vehicles/bicycles or on foot on a freelance or contract basis \citep{verlinde2016favor, archetti2021recent, FLECKENSTEIN2023499}. In Figure \ref{fig:stakeholders}, we present the connections among the stakeholders.

\subsubsection{Process Flows in SDD}
In SDD, customer orders are received continuously throughout distribution operations, which makes it necessary to develop a rapid decision-making system regarding request acceptance and parcel dispatch. When an order arrives, it triggers the process that matches the order with delivery logistics including delivery time, locations, vehicles, or drivers. If not matched, the request is either rejected or scheduled for delivery the next day. In the case of a successful match, the service provider proceeds to make dispatching decisions with a delivery performed by either drivers of the SDD company, specialized third-party logistics providers, or crowdsourced couriers. Subsequently, routing decisions are made. The process iterates until the end of the time horizon. Objectives might consider optimizing different performance metrics, such as maximizing expected request fulfillment, maximizing profit, minimizing vehicle operating costs, or minimizing $CO_2$ emissions as part of environment sustainability goals\footnote{https://fortune.com/2019/03/01/amazon-day-delivery-emissions/}. The process is summarized in Figure \ref{fig:SDDprocess}. 
Detailed decision specifications are provided below. 


We define \textit{Matching} as the agreement between the customer and the service provider on when and where to deliver. For each new order arrival, the service provider makes a feasibility check to determine if there is available delivery capacity to serve the customer. If the feasibility check yields a positive result, the customer is typically offered a selection of delivery deadlines, time windows, location options, and associated prices. Subsequently, the customer chooses a preferred delivery time and location if interested, and afterwards the customer request is assigned to a driver. In case the service provider cannot find feasible options to offer or if the customer is unsatisfied with all options, the request is rejected or scheduled for the next day. If matched successfully, both the customer and the driver are notified.


\textit{Dispatching and Vehicle routing} is determining the vehicle routes and schedules for fulfilling the available orders. Upon capturing and matching an order (see Decision 1 in Figure \ref{fig:SDDprocess}), the service provider faces dispatching decisions, wherein they decide whether to dispatch all or a subset of the currently available orders immediately or await additional orders, for potential consolidation purposes (see Decision 2 in Figure \ref{fig:SDDprocess}). Once a dispatching decision is made, the subsequent step involves the routing of vehicles to their designated destinations. Notably, there exists various delivery means, for example,  trucks, drones, hybrid combinations of vehicles and drones, crowdsourced drivers, and others, with the selection depending on the available resources. 

It is essential to distinguish SDD from conventional or next-day delivery scenarios. In SDD, orders arrive continuously, rendering it impractical to wait for the aggregation of all parcel arrivals before initiating deliveries. Consequently, vehicles are often required to undertake multiple routes to meet the time windows stipulated by customers or to adhere to the delivery deadline. Thus for each dispatch, not all available parcels need to be delivered immediately, even when time permits, as some can be dispatched later together with incoming parcels that share close delivery locations, thereby optimizing overall time and distance.

As illustrated in Figure \ref{fig:stakeholders}, there are multiple stakeholders in SDD routing problems, each with their own concerns and objectives. However, unlike previous surveys (e.g., \cite{WAMUTH2023801}, \cite{FLECKENSTEIN2023499}), we primarily focus on logistic operations that occur after the customer interaction, namely, the routing and distribution aspects (e.g., Decision 2 in Figure \ref{fig:SDDprocess}), as emphasized in the title of our study. The mixed-integer linear programs presented in the following sections pertain specifically to the decisions made by the service company, i.e.,  the service company decides which location options or time windows to offer to each customer. 

\subsubsection{Sources of uncertainty and handling measures}
SDD systems have to deal with various sources of uncertainty. In Figure \ref{fig:SDDprocess}, we show the stages within the SDD process where uncertainty arises, and we list the corresponding factors of uncertainty below the stages. 

When considering order arrivals, the common uncertainty factors are arrival time, demand quantity, order volume/weight, and delivery/pickup location. During the matching decision phase, there is uncertainty in customers' choices on delivery time slots or drivers' preferences on the requests they would like to serve especially when they are given flexibility to choose. 
When making dispatching and routing decisions, vehicle/courier availabilities are frequently dynamic and stochastic. Furthermore, during routing processes, external factors, such as traffic congestion, may induce variations in vehicle travel time. In scenarios where parcels are delivered to pickup stations, the time customers pick up their parcels is uncertain, which matters because of the limited capacity of pickup stations. 

As mentioned above, in this article we provide several deterministic mathematical models. These models can accommodate the above-mentioned uncertainties in different ways:
\begin{itemize}
    \item  By employing scenario sampling where a scenario represents a realization of uncertainty factors, we can solve a deterministic model for each scenario and then get final solutions by aggregation through e.g., consensus functions \citep{voccia2019same, li2022reinforcement}; 
    \item Deterministic models can be used to generate \emph{a priori} solutions under specific assumptions about the uncertain factors. These solutions may then be updated when new information arrives. Alternatively, one may opt for immediate implementation of the first action and recalculate a new \emph{a priori} policy in response to updated information, employing a dynamic rollout scheme that iterates through the process (\citet{klapp2018dynamic}).
\end{itemize}

In the papers we reviewed, a predominant focus lies in addressing dynamic and stochastic problems, often modeled as sequential decision processes. The components of such a decision process are: decision epochs; an information model describing the uncertainty over time; states describing the information available at each decision epoch, and decisions representing solutions of the decision model instance. At each decision epoch, a realization of the information model becomes known. To find an optimal solution is to find an optimal policy maximizing the expected sum of rewards. The modeling of a sequential decision process includes three steps: transferring observed uncertainty to an information model, defining the decision model, and modeling their interactions (see, e.g., \cite{soeffker2021stochastic}). 

To transfer uncertainty into an information model, suitable distributions must be adopted for each uncertainty factor. 
For instance, Poisson processes, Gaussian, and uniform distributions are commonly employed to model request arrival times (like e.g., \citet{ulmer2021restaurant, stroh2022tactical}), while uniform distributions are frequently assumed for request locations (\citet{ulmer2017delivery}). These distribution assumptions predominantly serve for generating scenarios in solution frameworks. However, uncertainty factors can be modeled with additional considerations to be more applicable to real-world scenarios. For example, one can define a time-varying order rate over the service day, following that customers are more likely to place orders near the cut-off time (\cite{banerjee2022fleet}), or time correlation can be introduced to travel time \citep{prokhorchuk2019stochastic, pureza2008waiting}.
Besides, distributions concerning vehicle or courier availability must describe rare events (\cite{soeffker2021stochastic}). A driver's availability may show correlations with realized and expected demand, while the distribution of vehicles may consider factors like weather or road conditions. Also, there are a few papers that integrate logit choice models into their solution frameworks to simulate customer behavior (\cite{prokhorchuk2019stochastic}).  

The solution approaches to derive policies from a sequential decision process depend on the exploitation of the information model. They include rolling horizon heuristics,
policy function approximation (PFA), cost function approximation (CFA), scenario sampling, rollout algorithm, value function approximation (VFA) or Q-learning. 
We refer to the survey by \cite{soeffker2021stochastic} for an exhaustive overview of information models and the advantages and disadvantages of the various approaches.


\subsection{Classification scheme describing characteristics of SDD problems}
\label{classification_scheme}

\begin{figure}[!h]  \centering\includegraphics[width=0.9\textwidth]{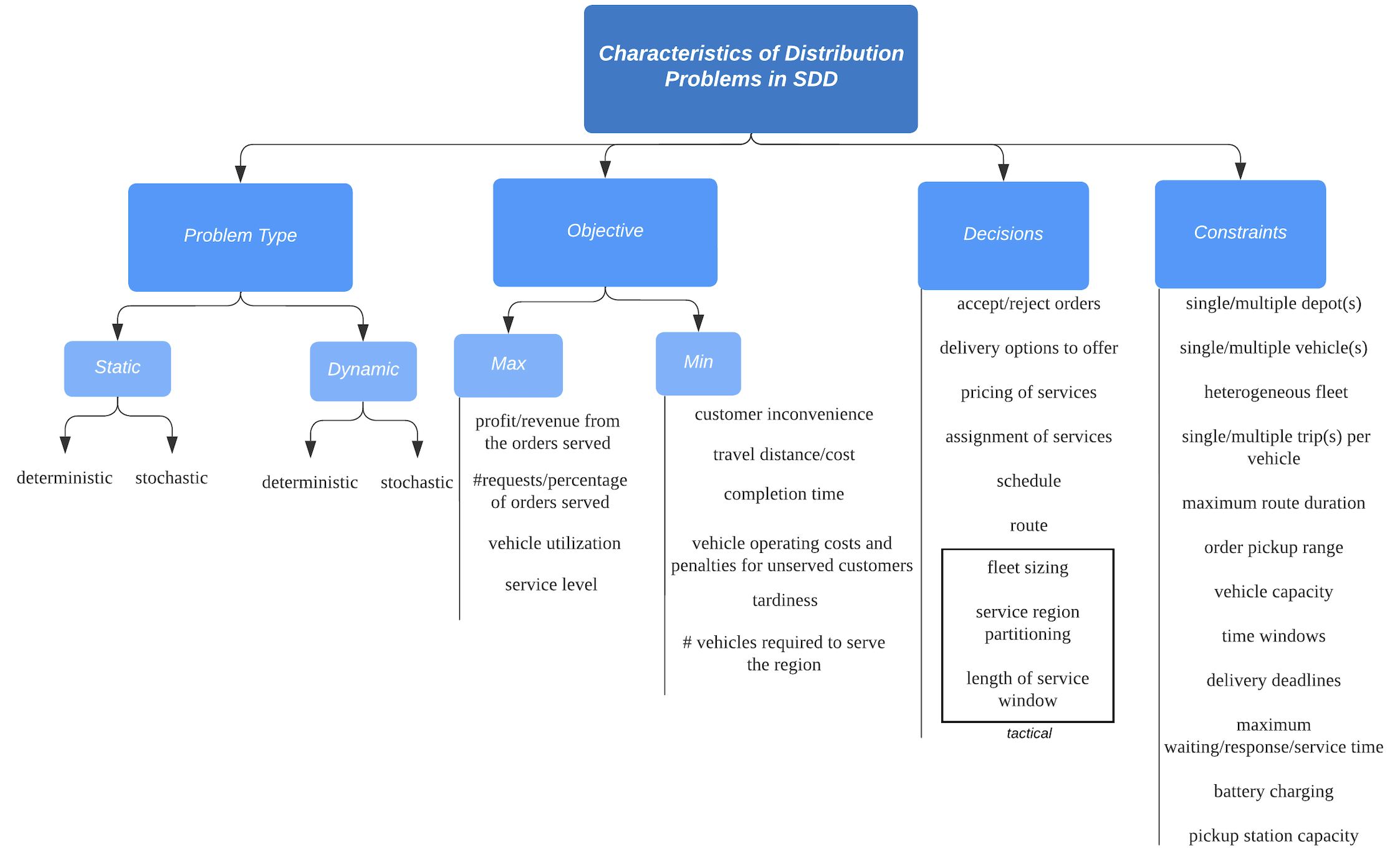}
    \caption{Classification scheme describing characteristics of SDD problems}
    \label{fig:SDDscheme}
\end{figure}

In Figure \ref{fig:SDDscheme} we provide a classification of SDD problem characteristics, namely, problem types, typical objective functions, decisions, and constraints. 

For problem types, we follow the classification defined by \cite{psaraftis2016dynamic}: 
\begin{itemize}
    \item Static and Deterministic: all inputs are known with certainty before determining solutions; there are no stochastic inputs.
    \item Static and Stochastic: solutions are still determined a priori with fixed parameters, but some elements follow probabilistic patterns.
    \item Dynamic and Deterministic: inputs are received dynamically, but no stochastic information is known about future inputs. 
    \item Dynamic and Stochastic: inputs are updated dynamically with stochastic information about future events.
\end{itemize}   

The papers we reviewed mostly handle operational challenges, including the optimization of customer service and financial goals. Specifically, the objectives commonly explored, such as the minimization of delays and the maximization of the number of requests served, are proxies for maximizing customer service; the objectives such as the maximization of profit/revenue and the minimization of costs represent the companies' financial needs. Furthermore, at the tactical level, it is common to focus on the efficient allocation of resources, such as minimizing the number of vehicles required to serve the region. 

The three decisions listed in the box in Figure \ref{fig:SDDscheme} are usually investigated at the tactical level and they are often regarded as fixed parameters when addressing operational problems. The remaining decisions shown in Figure \ref{fig:SDDscheme}, displayed sequentially to align with the general SDD process illustrated in Figure \ref{fig:SDDprocess}, are at the operational level.

At the tactical level, the company must determine fleet size or capacities distributed across different service zones (see Section \ref{Districting}). The realm of tactical decisions leads to a set of prevalent optimization questions, as delineated in recent works \citep{stroh2022tactical, banerjee2022fleet, wu2022service}: 
\begin{itemize}
    \item Fleet sizing - What is the optimal fleet size and composition?
    \item Order cutoff time - What is the latest feasible time for offering SDD services?
    \item Sizing of service area - What is the appropriate scale for defining a service area?
    \item Service network design - How to determine the shipment paths in both temporal and spatial dimensions?
\end{itemize}  
Typically, these models are constructed to capture the average system behavior in consideration of sampled or forecasted demand \citep{stroh2022tactical, banerjee2022fleet}. While the specific objectives of papers addressing these inquiries may differ based on their specific areas of focus, their overall goal falls into developing simplified models capable of facilitating efficient day-to-day operational activities at minimal costs.

Regarding constraints, unlike conventional delivery, vehicle capacity is not of utmost importance but delivery deadlines are, as SDD systems face low order volume and tight time to deliver (\cite{stroh2022tactical}). There are also context-specific but trending constraints, such as pickup station capacity when considering delivering to pickup stations, battery charging in the context of electrical vehicles or drones, or order pickup range when addressing customer preferences regarding the maximum distance they are willing to travel to pick up their orders from alternative locations. 


\section{Release dates and deadlines}
\label{RD_deadlines}
In traditional delivery settings, all orders are available at the depot before deliveries. On the contrary, in SDD, 
orders arrive throughout the day and must be served by the end of the day. \emph{Release dates}  is the term used to represent either the time at which the parcels associated with orders become available at the depot, or the order time (in case parcels are already available at the depot, which acts as a warehouse).
The products ordered online must be processed immediately by the distribution planners, because of tight delivery deadlines, hence the fulfillment operations often happen simultaneously with the distribution operations. 
Moreover, the tight delivery time complicates logistics management as consolidation opportunities are significantly reduced. 
This raises challenging questions: should the distribution planner wait before dispatching the vehicles for delivery to have more parcels available at the depot (and so to increase consolidation opportunities), or should the vehicles deliver the parcels that are already available in order not to miss deadlines? Also, if the distribution planner decides to wait, then for how long? 

To answer the above questions, a family of \emph{routing problems with release dates and/or deadlines} has emerged in the recent literature. These problems aim to model the fundamental tradeoff between time and cost in SDD. In the following, we first present a mathematical model for the basic setting with deterministic input data in Section \ref{formu_RD}, and then we discuss the most recent modeling and solution approaches that incorporate data uncertainty in Section \ref{lite:RD}. 
Note that these problems are introductory to the problems presented in the following sections. In this section, we present the works that do not have overlapping features with the classes of problems analyzed in the following sections. 

\subsection{Mathematical formulation}
\label{formu_RD}
 In this section, we focus on the deterministic setting and present the formulation for basic routing problems with release dates and/or deadlines. As mentioned in Section \ref{intro}, our goal is to analyze the emerging optimization problems in SDD and present mathematical formulations for the basic problem in each class. This formulation might be considered as a reference for researchers dealing with more sophisticated problems in the class and can be used as a base point, either for comparisons or possible generalizations. 
 For each class of problems analyzed, we present the 
 underlying deterministic counterpart for the single vehicle case.

In the class of routing problems with release dates, we specifically focus on the deterministic Orienteering Problem with release dates (OP-rd) introduced in \citet{li2022reinforcement}. {This is a basic problem in the class of routing problems with release dates whose setting and formulation, can be easily adapted to different variants and to more general cases. }The problem is described as follows.

Let us assume that we have a set of parcels $N$ with known release dates, $r_i\ge 0$, $i \in N$.  {Each parcel corresponds to a customer order and has to be delivered to the customer's place. The release date $r_i$ represents either the time at which the order is placed (in case the parcel is already at the depot), or the time at which the parcel becomes available at the depot (for an earlier placed order). A single vehicle is available to perform delivery operations. The vehicle makes multiple trips that all start and end at the depot. During each trip, the vehicle delivers a subset of the parcels available at the depot (thus, those associated with a release date not larger than the starting time of the trip).} 
We denote as $V$ the set of locations including the customers' locations and the depot, denoted as 0. Set $A$ is the set of arcs composed by any pair of nodes $i, j$ in $V$. The travel distance between two nodes is denoted as $d_{ij} \geq 0$, where $i, j \in V$.  {Distances satisfy the triangle inequality.} The goal of the problem is to maximize the number of customers that are served within the time horizon $[0,T_E]$, where $T_E$ represents a \emph{deadline} and can be interpreted as the end of the driver's shift.  {As mentioned above, the vehicle performs a set of trips $\cK$ to deliver the orders within the deadline, where a trip is a route that starts and ends at the depot and visits some customers in between.} The solution is composed of $|\cK|$ trips, where a trivial upper bound on the value of $|\cK|$ is given by the number of orders. A tighter upper bound is proposed in \cite{li2022reinforcement}. The trips $k\in \cK$ are numbered starting from 0. 
The OP-rd is NP-hard as it contains the OP as a special case in which all release dates are set to 0.  A mathematical formulation for the problem is proposed in \cite{li2022reinforcement}.
The problem is modeled as a mixed-integer linear program (MILP) using the following decision variables:
\begin{enumerate*}[(i)]
    \item $y_i^k$: binary variable equal to 1 if parcel $i$ is delivered in trip $k$, 0 otherwise;
    \item $x^k_{ij}$: binary variable equal to 1 if trip $k$ traverses arc $(i,j)$, 0 otherwise;
    \item $s_k$: continuous variable representing the starting time of trip $k$.
\end{enumerate*}
The formulation for the OP-rd is:
\small
\begin{subequations}
	\label{eq:base1}
\begin{equation}\label{eq:obj_ILP_base}
\hbox{max} \sum_{i \in N}\sum_{k \in \cK} y^k_i 
\end{equation}
\begin{align}
\label{eq:inout_ILP_base}
 \sum_{(i, j)\in A}x^k_{ij} &= \sum_{(j, i)\in  A}x^k_{ji}  = y_i^k && \quad i \in V, k \in \cK\\
\label{eq:subtour_ILP_base}
 \sum_{i, j\in S}x_{ij}^k & \leq \sum_{i \in S\setminus{\ell}}y_i^k 
&& \quad S \subseteq  N, |S|\geq 2, \ell \in S, k \in \cK \\
\label{eq:singleVisit_ILP_base}
\sum_{k \in \cK}y_i^k & \leq 1 &&  \quad i \in N \\
\label{eq:releaseDate_ILP_base}
s_k & \ge r_i y_i^k &&  \quad i \in N, k \in \cK \\
\label{eq:consecutiveTrips_ILP_base}
s_{k+1} & = s_k + \sum_{(i, j)\in A}d_{ij} x^k_{ij}  && \quad k \in \cK \\
\label{eq:finalTrip_ILP_base}
s_{|\cK|} & = T_E &&  \\
{x^k_{ij}, y^k_{i}} & \in \{0,1\} && i,j \in V, k \in \cK.
\end{align}
\end{subequations}
\normalsize
The objective function \eqref{eq:obj_ILP_base} maximizes the number of customers served. Equations \eqref{eq:inout_ILP_base} are in- and out-degree constraints, inequalities \eqref{eq:subtour_ILP_base} are the generalized subtour elimination constraints, and  \eqref{eq:singleVisit_ILP_base} state that each customer can be served at most once. Constraints \eqref{eq:releaseDate_ILP_base}--\eqref{eq:finalTrip_ILP_base} determine when each trip starts. Specifically, constraints \eqref{eq:releaseDate_ILP_base}
force a trip to start after the arrival of parcels delivered in the trip, constraints \eqref{eq:consecutiveTrips_ILP_base} ensure that a new trip starts when the previous trip finishes, and constraints \eqref{eq:finalTrip_ILP_base} state that the last trip ends at $T_E$. 
The equality in constraints \eqref{eq:consecutiveTrips_ILP_base} and  \eqref{eq:finalTrip_ILP_base} is due to the fact that we can shift the waiting time to the start of the distribution without affecting the structure and the solution value so that any waiting time between consecutive trips and the end of the shift is avoided (see \cite{archetti2018iterated} for more details).

Formulation (\ref{eq:base1}) can be adapted to more general or similar problems. For example, instead of maximizing the number of parcels served, the objective could be to minimize the completion time for serving all customers (see \cite{archetti2020dynamic}) by forcing each customer to be visited (changing inequality constraints (\ref{eq:singleVisit_ILP_base}) to equality) and removing constraints (\ref{eq:finalTrip_ILP_base}). We can also handle the problem with multiple capacitated vehicles by adding capacity constraints and applying the constraints for each vehicle. Problems with additional side constraints, like time windows for the deliveries, can be easily integrated as well.

Most of the contributions mentioned in the following subsection deal with the stochastic and dynamic setting, in which the uncertainty is mostly related to release dates. The problem is then typically modeled as a Markov Decision Process (MDP)  \citep{ulmer2019same,van2019delivery,archetti2020dynamic,ulmer2020dynamic,li2022reinforcement}.
When dealing with dynamic and stochastic problems, Formulation \eqref{eq:base1} can be used in multiple ways: for competitive analysis (i.e., evaluating the solution with perfect information), or for finding a solution in each decision epoch using an approximation 
of future information. Both ways are applied in \cite{li2022reinforcement}, so we refer the reader to this paper for a reference about using a basic deterministic model for the solution of a stochastic and dynamic problem.


\subsection{Literature review on SDD routing problems with release dates and deadlines}\label{lite:RD}

Some early contributions exist in studying routing problems where parcels to be delivered arrive dynamically during the distribution process. Though the authors did not use the term ``release date'', the setting corresponds to the one introduced above. In accordance with the concepts introduced in Section \ref{classification_scheme},  we summarize the papers to give an overview and provide features of Model (\ref{eq:base1}) as a point of reference in Table \ref{table:literature_RD}. The citations in Table \ref{table:literature_RD} are listed in the same order that appeared in the text. We observe that most papers deal with dynamic and stochastic VRP or TSP aiming to maximize the number of requests served or minimize travel distance/cost. Besides, the common solution approach employed is ADP.

\begin{table}[!h]
\caption{Summary of papers on SDD routing problems with release dates and deadlines. The abbreviations: adaptive large neighborhood search (ALNS), tabu search (TS), variable neighborhood search (VNS) multiple depots (MD), single-vehicle (SV), multiple trips per vehicle (MT), vehicle capacity (VC), time windows (TW), maximum waiting/response/service time (WT), delivery deadlines (DD, which denotes either the end of drivers' shift or the latest allowable time for delivering to customers).}

\label{table:literature_RD}
\scalebox{0.61}{\color{black}
\begin{tabular}{p{0.2\textwidth}p{0.1\textwidth}p{0.23\textwidth}p{0.2\textwidth}lp{0.2\textwidth}lp{0.14\textwidth}}\hline
\textbf{Literature}     & \textbf{Problem}    & \textbf{Objective}                                                                          & \textbf{Uncertainty}                                                                       & \textbf{Uncertainty Nature} & \textbf{Decision}              & \textbf{Constraint} & \textbf{Solution Approach}            \\ \hline
\cite{azi2012dynamic}         & VRP                 & max(profit collected from the orders served)                                            & requests(arrival time, location)                                                          & dynamic \& stochastic                        & accept/reject; route           & RD; TW; MT; DD      & ALNS                         \\
\cite{ferrucci2013pro}   & VRP                 & min(customer inconvenience)                                                                 & requests(arrival time, location)                                                          & dynamic \& stochastic                        & route                          & TW; DD; WT          & TS; VNS                           \\
\cite{ferrucci2016pro}   & VRP                 & min(customer inconvenience)                                                                 & requests(arrival time, location)                                                          & dynamic \& stochastic                        & route                          & TW; DD; WT          & ML; TS                  \\
\cite{cattaruzza2016multi}  & VRP                 & min(travel distance)                                                                        &         no                                                                                   & static \& deterministic                         & route; assignment of services             & TW; VC; DD; MT     & memetic algorithm          \\
\cite{yao2019robust}         & TSP                 & min(travel cost)                                                                            & requests(location)                                                                            & dynamic \& stochastic                        & route                         & DD; MT;SV           & MILP         \\
\cite{van2019delivery} & delivery dispatching                 & min(travel distance)                                                                        & requests(arrival time, demand quantity, order volume, demand location, dispatch window) & dynamic \& stochastic         & schedule; route                & TW; VC; DD; MT     & ADP; continuous approximation                          \\
\cite{voccia2019same}       & VRP                 & max(\#requests served)                                                                      & requests(arrival time)                                                                    & dynamic \& stochastic         & schedule; route               & TW; DD; MT         & scenario sampling   \\
\cite{cote2023branch}       & VRP                 & min(\#requests rejected) and min(routing cost)                                                                      & requests(arrival time)                                                                    & dynamic \& stochastic         & schedule; route; preemptive return               & TW; DD; MT         & branch-and-regret, scenario sampling; local search; ALNS; Regret-k heuristic \\
\cite{archetti2020dynamic}    & TSP                 & min(completion time for serving all customers)                                              & requests(arrival time)                                                                    & dynamic \& stochastic         & schedule; route                & MT; SV              & local search                 \\
\cite{li2022reinforcement}         & OP                  & max(\#requests served)                                                                      & requests(arrival time)                                                                    & dynamic \& stochastic         & schedule; route               & DD; MT;SV           & RL; continuous approximation                           \\
\cite{klapp2018dynamic}       & TSP                 & min(vehicle operating costs and penalties for the unserved)                                 & requests(arrival time)                                                                    & dynamic \& stochastic         & schedule; route                & DD; MT;SV           & ADP; local search                          \\
\cite{ulmer2017delivery}            & VRP                 & min(sum of delay over all customers)                                                        & requests(arrival time, location)                                                          & dynamic \& stochastic         & route                         & TW; DD; MD         & cheapest insertion           \\
\cite{ulmer2019preemptive}      & VRP                 & max(\#requests served)                                                                      & requests(arrival time, location)                                                          & dynamic \& stochastic         & preemptive return; route            & DD; SV; MT          & ADP                          \\
\cite{klapp2020request}      & TSP                 & min(vehicle operating costs and penalties for the unserved)                                 & requests(arrival time)                                                                    & dynamic \& stochastic         & accept/reject; schedule; route & DD; SV; MT          & ADP        \\
\cite{leung2022community}     & VRP & max(usage efficiency of   vehicles) and min(the waiting time of orders) & requests(arrival time, location, weight)                                                  & dynamic \& stochastic         & schedule; route                & DD; VC; WT          & community logistics strategy\\\hline
Model (\ref{eq:base1})  & OP & max(\#requests served) & no                                                 & static \& deterministic         & schedule; route                & DD; SV; WT          & MILP
 \\\hline
\end{tabular}
}
\end{table}

{\cite{azi2012dynamic} study a problem arising in the distribution of online grocery orders. Orders arrive dynamically during the day and the provider has to decide, each time an order is placed, whether to accept it or not. Each order is associated with a profit and a fleet of vehicles is available to perform delivery operations. Each vehicle can perform multiple trips of a limited duration. The goal is to maximize the profit collected from the orders served. The authors propose an Adaptive Large Neighborhood Search (ALNS) embedding an acceptance rule that takes into account potential future requests. The algorithm is compared against a myopic rule and the results show that it provides a competitive advantage.}

{\cite{ferrucci2013pro} consider a daily distribution process evolving in real-time, such as the delivery of newspapers. Customer requests arrive dynamically during the day. Each request is associated with a soft time window, and a high penalty is imposed if service begins after it. 
The problem is to decide at each specified time point the tour plan that assigns the set of unserved requests to vehicles, with the goal of minimizing the total customer inconvenience that is captured by a function of request response times. The authors exploit past request information to anticipate future requests. They integrate these anticipated requests with known requests to decide which customers should be served and then use a Tabu Search (TS) with a variable neighborhood search to generate the tour plan. 
The authors prove the efficiency of their approach by comparing it with a deterministic approach.}

{\cite{ferrucci2013pro} assume request arrivals on different days follow the same spatial and temporal structure, thus they use a single profile of past requests to anticipate future request arrivals. However, \cite{ferrucci2016pro} observe that there are days that seem to follow a different structure, so they extend the single profile approach by applying multiple profiles for the same problem setting studied in \cite{ferrucci2013pro}. They identify $k$ profiles by applying a $k$-means algorithm (\cite{macqueen1967classification}). 
During the day, request arrivals are analyzed in real time for mapping the request arrival structure to a profile. The computational results compared to the previous single profile approach prove the superiority of the new approach.}

\citet{yao2019robust} consider a setting in which a single uncapacitated vehicle is dispatched from a local distribution center, with only a subset of parcels' locations being known in advance. Additional parcels' locations become known while the vehicle is performing deliveries and preemptive returns are allowed. Preemptive return means that the vehicle might return to the depot before having delivered all parcels onboard, in case this is convenient.   
The problem is modeled as a precedence-constrained asymmetric traveling salesman problem with demand uncertainties. The uncertainty is addressed with robust optimization and a MILP formulation is proposed. 

{To the best of our knowledge, the first paper using the term `routing problems with release dates' is \cite{cattaruzza2016multi}}. The authors introduce the problem in the context of a two-level distribution system in which, on the first level, parcels are transported to the local distribution center by trucks. Release dates represent the time at which the parcels reach the distribution center. It is assumed that release dates are known in advance, i.e., they are deterministic. In the second level, a  fleet of vans is used to distribute goods to the final customers.
The parcels continue to arrive during the day at the distribution center and 
each parcel has to be delivered within a pre-defined time window. 
 The goal is to determine the distribution plan for the vans so that all parcels are delivered within the requested time windows and the routing cost is minimized. The authors propose a memetic algorithm and introduce a set of instances to serve as a benchmark for future works.


While the study by \cite{cattaruzza2016multi} addresses a problem within a deterministic setting, more recent literature has focused on variations using stochastic and dynamic input data, often modeled through MDP. We now review the most prominent contributions among these studies. 

\cite{van2019delivery} present a variant of the {multi-vehicle} delivery dispatching problem. In their setting, the number of orders arriving, their volume, destinations, and dispatch windows are uncertain. Every order must be shipped within its dispatch window. Dispatch operations are performed using a fleet of primary and secondary vehicles, where secondary ones are more expensive, but it is assumed that all orders can be dispatched within the planning horizon.  
The goal is to minimize the expected dispatching costs. 
The authors model the problem as an MDP. To handle  large instances, 
they propose an approximate dynamic programming algorithm (with the approximation of routing cost) where the downstream costs are estimated through a linear value function approximation.

\cite{voccia2019same} consider a multi-vehicle  delivery problem where requests  arrive dynamically. The customer locations are known in advance, but the release dates are subject to uncertainty. It is assumed that the arrival rate and the distribution of release dates for each order are given. Moreover, for each order, we are given a (deterministic) deadline or a time window during which the delivery
must occur. Deliveries are performed by a given fleet of uncapacitated vehicles, and no preemptive returns are allowed. 
The goal is to maximize the expected number of requests that are served on time before the end of the planning horizon. The authors first provide an MDP model for the problem. By incorporating information about future requests into routing decisions, they then derive 
an analytical result that identifies when it is beneficial
for vehicles to wait at the depot. 
They also present a heuristic solution approach based
on scenario sampling combined with a consensus function. 

Similarly, \cite{cote2023branch} study a multi-vehicle problem with stochastic and dynamic requests and time windows. Their primary objective is to minimize the number of rejected requests, while the secondary objective is to minimize the total routing cost. Besides, they accept preemptive returns. To address the problem, they propose a novel branch-and-regret (B\&R) heuristic, which incorporates scenario sampling techniques, consensus functions, and a branching scheme. Through testing on the same instances as the ones utilized by \cite{voccia2019same}, the authors show that their B\&R performs much faster, while also achieving a higher number of served requests.

\cite{archetti2020dynamic} present a single-vehicle problem with stochastic and dynamic release dates and
model the problem as an MDP. 
The objective is to minimize the total time for serving all customers. They propose a reoptimization approach that heuristically reduces the number of decision epochs. The optimization problem at each decision epoch is solved in two ways: through a deterministic model that takes as input the point estimation of release dates corresponding to the expected value; and through a stochastic optimization model that considers the entire stochastic information on the release dates. At each reoptimization epoch, they derive the solution of the two models through an iterated local search. The computational results show that the stochastic model performs better than the deterministic one but is more time-consuming. The deterministic counterpart of their problem has been studied earlier by \citet{archetti2018iterated}.

\citet{li2022reinforcement} introduce the orienteering problem with stochastic and dynamic release dates, whose deterministic counterpart corresponds to the problem presented in Section \ref{formu_RD}. The problem can be seen as the single vehicle version of the problem studied in \citet{voccia2019same}. In the setting analyzed in \citet{li2022reinforcement}, the release dates
are stochastic and follow a known distribution. As in \citet{voccia2019same}, no preemptive returns are allowed, and the
problem is modeled as an MDP with the goal of maximizing the expected number of requests served before the end of the planning horizon. 
The authors propose an exact branch-and-cut method for the solution of the optimization problem at each  decision epoch (as opposed to heuristic methods used by \citet{voccia2019same}). 
Specifically, they propose an approximation of future states and derive decision policies through value function and policy function approximation. Branch-and-cut exact methods are used to evaluate individual scenarios before applying a consensus function (in the policy function approximation), and a two-stage stochastic MILP model is proposed to simultaneously take multiple sampled scenarios into account (in the value function approximation). 
The authors empirically demonstrate the advantage of embedding exact methods inside a reinforcement learning
framework, to improve the sequential decision-making process.

As pointed out in the works of \citet{voccia2019same} and \citet{li2022reinforcement},  in SDD not all orders may be served by the end of the day. To deal with unserved orders, \cite{klapp2018dynamic} impose penalties and  minimize the (expected) vehicle travel cost together with penalties for unserved requests. 
The problem is modeled as an MDP, where dispatching decisions are made at fixed intervals and no preemptive returns are allowed. As in \citet{li2022reinforcement}, a single vehicle is considered and 
all requests have the same deadline, i.e., the end of the planning horizon.  In addition, it is assumed that all customers lie on a line segment and that the depot is located at one end of this segment.
The authors first propose 
an integer programming model 
and a local search heuristic to solve the deterministic variant of the problem where all the release dates 
are known. 
Then they use this approach to derive \emph{a priori} solutions for the stochastic variant 
with predetermined routes.
They finally derive heuristic solutions and lower bounds for the more general case with stochastic release dates. 


\cite{ulmer2017delivery}  studies a dynamic multi-vehicle multi-depot routing problem with stochastic release dates and models it as an MDP.
The problem setting includes  \emph{delivery deadlines}, which are defined as time windows starting from the uncertain release date (the time when the customer places the order) and whose duration is the delivery time promised to the customer. 
Contrary to \citet{voccia2019same}, where hard deadlines are imposed (i.e., the customer is considered ``served'' only if the delivery happens within the given time window), \citet{ulmer2017delivery} allows for \emph{soft deadlines} (i.e., in case the order is not delivered before the deadline, the customer experiences a delay).  
The goal is to minimize the expected average delay per customer. The author performs a simulation study on the Iowa City area and concludes that narrow delivery deadlines cause substantially higher delivery costs.

\cite{ulmer2019preemptive} investigate a similar setting as in \citet{li2022reinforcement}, with a single vehicle performing the delivery, with stochastic release dates and with the aim of maximizing the expected number of served customers within the given planning horizon. However, as in \citet{yao2019robust}, the authors allow preemptive returns to the depot. They use approximate dynamic programming as a solution method. 
Through an extensive computational study, they show that preemptive returns are most effective when the returns occur late enough in the planning horizon, while they are much less effective in other cases.   

\cite{klapp2020request} study a problem in which customers place orders over time, and the dispatcher immediately decides whether to accept the order (and thus, guarantee that the parcel will be delivered before the end of the day) or to deny the service. 
A single vehicle is used to perform deliveries, no preemptive returns are allowed, and a penalty is imposed for each denied request. The goal is to determine the distribution plan that minimizes the expected sum of vehicle travel costs and penalties for denials.
Vehicle routes are updated
dynamically and serve each accepted delivery request no later than the end of the service day.

In the models presented so far, consolidation opportunities  are implicitly exploited, while considering certain objectives (like the maximum number of deliveries, minimum delivery cost, and penalization of unserved deliveries). 
A more explicit objective is proposed by \cite{leung2022community}. The authors consider a delivery problem combining multiple vehicles with pickup stations. The problem is to decide when each vehicle departs from the depot and which orders are delivered by each vehicle. The goal is to maximize the usage efficiency of vehicles and minimize the waiting time of orders before delivery. The authors propose an approach that clusters packages into communities where each community is served by one vehicle. Then decisions associated with each community are made in order to better utilize the available vehicle capacity according to different criteria: (i) time - temporal delivery postponement, determine the most suitable shift for each customer request; (ii) space - spatial community adjustment, adjust the size of the geographical area served by every vehicle.  The computational results show that the spatial approach works better in terms of  order fulfillment rate and the value of postponement might not be enough to compensate for its cost.

\section{Time window options for attended home delivery}
\label{time_window}

Many SDD companies offer customers the flexibility to decide when, where, and how their parcels should be delivered, so there will be fewer failed deliveries. Meanwhile, SDD companies face complex decision-making challenges, such as which time window options to offer, how to price different time window options, where to deliver when there are multiple possible delivery locations, and how to coordinate different transportation methods. In this section, we narrow our focus on time window options and pricing. Because of the importance of the temporal and pricing aspects in SDD, we provide a mathematical formulation of it, which warrants a separate and dedicated section away from the next one that discusses location and transportation options. 

In a recent UK consumer home delivery report, 1000 respondents were questioned on what they thought about the length of the delivery time slot\footnote{https://nshift.com/resources/delivery-time-slots-customer-satisfaction}. Over 40\% of the respondents stated they were willing to wait for a maximum of 2 hours, while less than 10\% would wait for more than 4 hours. This underlines the pressure from online customers on service quality: not only do they require fast deliveries (as mentioned in the former section), but also narrow time windows for delivery. 

To 
avoid or reduce delivery failures, many researchers recently investigated the option of offering multiple time windows to customers in attended home delivery (AHD) services.  This time window option is often applied in online grocery ordering. For each request, the customer and the grocery retailer must mutually agree on a time window for delivery. The process imposes computational challenges as the proposal of time window slots should be as fast as possible (within milliseconds) to ensure a smooth booking process in SDD services.  Besides, the combination of SDD and narrow deadlines is economically challenging for service providers (\cite{UKR}). When the deadline is tight, vehicles often only deliver a few orders per route. As a consequence, the revenue is low. Thus, another research trend is emerging that focuses on studying \textit{delivery time options with pricing}, i.e., customers are offered different delivery times associated with different costs (the tighter the delivery time, the higher the cost). In the following, we will first present a mathematical model for the problem with delivery deadline options in a deterministic environment in Section \ref{formu_pricing}, then we will discuss recent works dealing with the challenges in deciding time slots or pricing options in Section \ref{pricing}.

\subsection{Mathematical formulation}
\label{formu_pricing}

As done in the former section, we now present the mathematical formulation of a basic problem arising in the case of delivery deadline options for each order. This formulation might be used as a benchmark or a starting point to deal with more complicated problems in this class. 

The formulation is built on the problem setting introduced by \cite{ulmer2020dynamic} but considering a single uncapacitated vehicle and a deterministic environment, i.e., the two unknown factors in \cite{ulmer2020dynamic} -  customer release dates and customer choice behavior based on willingness-to-pay (WTP) - are assumed to be known here. We describe the problem as follows.

{A set of customers $N$ order goods during the day, and we denote the order time of each customer $i$ as $r_i$ (which might be seen as the release date introduced in the former section). $V$ is the set of locations including the customers' locations $N$ and depot $0$, i.e., $V=N\cup 0$. A set of same-day delivery deadline options $\Delta$ is available. The corresponding deadline of each option is denoted with $D_\delta$ ($\delta \in \Delta$). For each option $\delta \in \Delta$, the service provider knows the highest amount a customer $i$ is willing to pay (also known as the WTP), denoted with $U_{\delta}^i$
 and has to determine a customer-specific price $p_{\delta}^i$.
If $p_{\delta}^i$ is smaller than $U_{\delta}^i$, the option $\delta$ is feasible for customer $i$.  If there are no feasible options, the customer will be moved to next-day delivery. 
Note that, as we are in a deterministic setting, for each option $\delta$ that is feasible for customer $i$, i.e., $p_{\delta}^i \leq U_{\delta}^i$,  $p_\delta^i$ should always be equal to $U_\delta^i$ in order to maximize the revenue. In this way, we can omit decision variable $p_\delta^i$ as the offered price will be set according to the delivery option which maximizes the revenue of the provider.
Thus, the provider selects the delivery option to be offered to each customer and builds the delivery routes for the vehicle, such that the revenue coming from the selected delivery options is maximized. Customers who are served on the next day are not considered in the delivery plan of the day. The vehicle makes maximally $|\cK|$ trips to serve the customers. 
All trips start from depot $0$, serve a  set of customers and go back to the depot. 
As in Formulation (\ref{eq:base1}), a solution consists of up to $|\cK|$ trips performed by the vehicle. The goal is to maximize the revenue from the selected delivery options within  $T=[0, T_E]$. 
}

{The following decision variables are used:
\begin{enumerate*}[(i)]
    \item $y_i$: binary variable equal to 1 if customer $i$ is served in any SDD deadline option, 0 otherwise (i.e., $i$ is served with next-day delivery);
    \item $y_{i\delta}^k$: binary variable equal to 1 if customer $i$ is served in trip $k$ with deadline option $\delta$, 0 otherwise;
    \item $x_{ij}^k$: binary variable equal to 1 in case arc $(i, j)$ is traversed in trip $k$, $i, j \in V$;
    \item $a_i$: continuous variable representing the arrival time at customer $i$; 
    \item $s_k$: continuous variable representing the starting time of trip $k$;
    \item $e_k$: continuous variable representing the ending time of trip $k$.
\end{enumerate*}
Furthermore, we introduce auxiliary variables $y_{0\delta}^k$ with $U^0_{\delta}=0$, for each $\delta \in \Delta$. 
}
The problem formulation is then:
\small
\begin{subequations}
	\label{eq:base2}
\begin{equation}\label{eq:obj_ILP2_base}
\hbox{max} \sum_{k \in \cK}\sum_{\delta \in \Delta}\sum_{i \in N} U^i_{\delta}y_{i\delta}^k
\end{equation}
\begin{align}
\label{eq:each_ILP2_most_once0}
\sum_{k \in \cK}\sum_{\delta \in \Delta} y_{i\delta}^k &= y_i &&\quad i \in N\\
\label{eq:each_ILP2_most_once1}
\sum_{\delta \in \Delta} y_{0\delta}^k &\leq 1 &&\quad k \in \cK\\
\label{eq:inout_ILP2_base}
 \sum_{(i, j)\in A}x^k_{ij} &= \sum_{(j, i)\in  A}x^k_{ji}  = \sum_{\delta \in \Delta}y_{i\delta}^k && \quad i \in V, k \in \cK\\
\label{eq:order0_ILP2}
    a_i & \geq r_i && \quad i \in N\\
\label{eq:order1_ILP2}
    a_i & \leq r_i + \sum_{\delta \in \Delta}\sum_{k \in \cK} D_{\delta}y_{i\delta}^k+M(1-y_i) && \quad i \in N\\
\label{eq:arrival_ILP2}
    a_j & \geq a_i +d_{ij}\sum_{k \in \cK}x_{ij}^k-M(1-\sum_{k \in \cK}x_{ij}^k) && \quad i,j \in N, i \neq j\\
\label{eq:depart_ILP2}
    a_j & \geq s_k + d_{0j}x_{0j}^k -M(1-x_{0j}^k) && \quad  j\in N, k \in \cK\\
\label{eq:as_ILP2}
    a_i & \leq e_k-d_{i0}x_{i0}^k+M(1-x_{i0}^k) && \quad  i\in N, k \in \cK\\
\label{eq:end_ILP2}
    e_k & = s_k +\sum_{(i, j) \in A}d_{ij}x_{ij}^k && \quad  k \in \cK\\
\label{eq:orderT_ILP2}
    s_k & \ge r_i \sum_{\delta \in \Delta}y_{i\delta}^k &&  \quad i \in N, k \in \cK \\
\label{eq:next_ILP2_base}
    s_{k+1} & \geq s_k + \sum_{(i, j) \in A}d_{ij}x_{ij}^k && \quad k \in \cK\\
\label{eq:e_ILP2}
    s_{|\cK|} & \leq T_E && \\
    x^k_{ij}, y^k_{i\delta} , y_i & \in \{0,1\}&& \quad i,j \in V, \delta \in \Delta, k \in \cK\\
    s_k, a_i, e_k & \geq 0 && \quad i \in V, k \in \cK 
\end{align}
\normalsize

The objective function (\ref{eq:obj_ILP2_base}) maximizes the revenue. Constraints (\ref{eq:each_ILP2_most_once0})  guarantee each customer is delivered at most once with one of the SDD options. 
Constraints (\ref{eq:each_ILP2_most_once1})  limit the depot to be visited at most once in each trip. Constraints (\ref{eq:inout_ILP2_base}) are degree constraints. Constraints (\ref{eq:order0_ILP2}) and (\ref{eq:order1_ILP2}) ensure that the delivery time of a parcel is after its order time but before the promised deadline. Constraints (\ref{eq:arrival_ILP2}) set delivery times of customers visited consecutively in a trip.  
Constraints (\ref{eq:depart_ILP2}) and  (\ref{eq:as_ILP2}) fix the starting and ending time of each trip with respect to the visiting time of the first and last customer in the trip. 
Constraints (\ref{eq:end_ILP2}) calculate the ending time of each trip on the basis of the starting time and the duration of the trip. Constraints (\ref{eq:orderT_ILP2}) ensure a trip can only start after the order time of the requests to be served. 
Constraints (\ref{eq:next_ILP2_base}) state that a new trip starts after the end of the preceding trip. Constraints (\ref{eq:e_ILP2}) ensure that the ending time is within the time horizon. Note that constraints (\ref{eq:arrival_ILP2}) are sufficient to prevent subtours. However, generalized subtour elimination constraints 
(\ref{eq:subtour_ILP2_base})
\small\begin{align}
 \sum_{i, j\in S}x_{ij}^k \leq \sum_{i \in S\setminus{\ell}}\sum_{\delta \in \Delta}y_{i\delta}^k, &&S \subseteq  N, |S|\geq 2, \ell \in S, k \in \cK &&\label{eq:subtour_ILP2_base}
 \end{align}
\end{subequations}\normalsize
can be introduced to strengthen the formulation by being added to the model directly when $|S|$ is small or dynamically when $|S|$ is large.

This formulation can be adapted to handle the case of multiple capacitated vehicles by applying the constraints above for each vehicle and adding capacity constraints.

The problem considered above is deterministic. The settings in the contributions mentioned in the following subsection are mostly stochastic and dynamic, where the release dates or customers' preferences are uncertain. The problem is then often modeled as an MDP \citep{prokhorchuk2019stochastic, ulmer2020dynamic, KLEIN2022}. When considering the uncertain arrival time of customer orders and the unknown customers' individual preferences (WTP), the Formulation (\ref{eq:base2}) can be used to evaluate the solution with perfect information for competitive analysis or it can be embedded into the stochastic and dynamic approach. 
When it comes to dynamic arrivals, 
the pricing rule should consider both the current revenue and the impact on the fleet's flexibility to serve future requests.

\subsection{Literature review on delivery time slots with pricing}
\label{pricing}

\begin{table}[!h]
\caption{Summary of papers on time window options for attended home delivery. The abbreviations: multiple trips per vehicle (MT), time windows (TW), single-vehicle (SV), and delivery deadlines (DD).}
\label{table:literature_timeW}
\scalebox{0.67}{\color{black}
\begin{tabular}{p{0.2\textwidth}llp{0.2\textwidth}lp{0.2\textwidth}lp{0.13\textwidth}}\hline
\textbf{Literature}     & \textbf{Problem}    & \textbf{Objective}                                                                          & \textbf{Uncertainty}                                                                       & \textbf{Uncertainty Nature} & \textbf{Decision}              & \textbf{Constraint} & \textbf{Solution Approach}            \\ \hline
\cite{prokhorchuk2019stochastic}   & VRP              & max(revenue)                               & requests(arrival time), travel time, customer choice & dynamic \& stochastic         & pricing of services; route                            & TW; DD; MT          & ADP; cheapest insertion \\
\cite{ulmer2020dynamic}                & VRP              & max(revenue)                               & requests(arrival time), customer choice               & dynamic \& stochastic         & delivery options to offer; pricing of services; route & DD; MT          & ADP                     \\
\cite{KLEIN2022} & VRP              & max(profits)                               & requests(arrival time), customer choice               &  dynamic \& stochastic         & delivery options to offer; pricing of services; route & TW; DD; MT          & ADP   \\\hline
Model (\ref{eq:base2})  & OP & max(revenue) & no                                                 & static \& deterministic         & delivery options to offer; route                & DD; SV; MT          & MILP                 \\ \hline
\end{tabular}}
\end{table}

In many AHD applications, once a customer initiates a delivery request, the provider {not only} derives anticipatory decisions on which delivery time slots should be offered to the customer, {but it may also decide to charge the customers for some of the offered time slots. Indeed, as shown by \cite{ulmer2017delivery}, delivering orders within a short time is costly. }
{The author quantifies the impact of delivery deadlines on delivery costs and finds that a deadline of 60 minutes may double the number of vehicles and working time needed to perform the service compared to a delivery service with a 3- or 4-hour deadline. However, not many users are willing to pay much for {instant deliveries.} }Based on a survey of over 4700 respondents in China, Germany, and the United States (\cite{joerss2016customer}), about 25\% of consumers (especially young people) are willing to pay premiums for the privilege of SDD or instant delivery. Nevertheless, only 2\% said they would pay sufficiently more to make the instant delivery viable. Thus, careful pricing management is necessary to meet the delivery costs while guaranteeing a certain level of customer satisfaction. Similarly, following the concepts introduced in Section 2.2,
we summarize the papers to give an overview and provide features of Model (\ref{eq:base2}) as a point of reference in Table \ref{table:literature_timeW}. We see that, besides the request arrival time, customer choice is a prevalent uncertainty factor considered as these problems involve direct customer interactions. 

\cite{prokhorchuk2019stochastic} investigate the problem of determining the price of different deadline options, by combining the optimization of pricing and routing. They consider uncertainty about travel times, which means the delivery may not arrive on time. They model the problem as an MDP. Pricing decisions are based on solving an optimization problem involving the customer choice model, probability of arriving on time, and penalties for late deliveries. Routing decisions are made via the cheapest insertion heuristic and then the best route is selected based on the expected profit.  {The results show that their algorithm outperforms {several} deterministic benchmark approaches.}

\cite{ulmer2020dynamic}  studies how to price different deadline options for customers arriving dynamically, similarly as in \cite{prokhorchuk2019stochastic}. They show how dynamic pricing can considerably increase companies' revenues and the number of requests served within the time horizon. {Unlike \cite{prokhorchuk2019stochastic} who consider uncertain travel time,  in \cite{ulmer2020dynamic}   travel times are deterministic and the focus is on the uncertainty of customer choice behavior, which is modeled through a distribution function of the customer's WTP for each delivery option.}  The authors first model the problem as an MDP and then propose an anticipatory pricing and routing policy that incentivizes customers to select delivery deadline options that are most convenient for the routing plan. The pricing policy considers both the current revenue and the opportunity cost. The latter is defined as the difference in future revenue between the case in which the customer accepts an option and when she does not. To consider the uncertain WTP, which is correlated with each customer and each deadline option, the state-dependent opportunity cost per customer and option is calculated. Then, an approximation through state space aggregation is implemented to deal with the curse of dimensionality.
 
Similarly to \cite{ulmer2020dynamic}, \cite{KLEIN2022}  optimize the delivery time slots and prices offered to each incoming customer request and assume deterministic travel time.   {The main difference between the two studies lies in the underlying solution approach.} {Unlike \cite{prokhorchuk2019stochastic} and \cite{ulmer2020dynamic}, whose methodology relies on learning-based value function approximation, \cite{KLEIN2022}  use scenario sampling (see \cite{bent2004scenario}) to derive an online tour plan and a choice-based pricing optimization problem with the discrete price points to determine options pricing.} The computational study shows that their approach outperforms a myopic benchmark approach.

\section{Delivery options}
\label{delivery}

{In this section, we analyze routing problems with delivery options arising in SDD services. The term \textit{delivery option} is used to indicate both the case where multiple options are offered to customers as to where the parcel should be delivered, and the case where different options are available concerning how the parcels are delivered (different transportation means, or different cohorts of drivers).}

E-retailers and shipping companies often provide different options to customers concerning where the parcel should be delivered (office, home, lockers, or car trunk\footnote{https://money.cnn.com/2018/04/24/technology/amazon-key-in-car-delivery-review/index.html}). 
Some shared delivery locations are, for instance, pickup stations (parcel locker) and shops (\cite{savelsbergh201650th}). They work as secured and unattended customer locations to which couriers deliver packages for later pickup by customers. 
In this way, customers have more flexibility to pick up their parcels rather than waiting at home. Moreover, drivers can deliver a bunch of parcels to a pickup station instead of serving customers individually, thus reducing operating costs {and pollution.}  Furthermore, failed deliveries are reduced, and this makes the constrained capacity in parcel-delivery networks less {tight}
(\cite{omnichannelDelivery}). 

{As for transportation means, autonomous vehicles can be considered as a novel technology to face the increasing shortage of drivers. In addition, crowdshipping is an emerging delivery strategy where ordinary people transport parcels instead of companies' drivers, against a compensation. The main advantage of crowdshipping is that it allows to dynamically adjust the delivery capacity according to demand fluctuations, without incurring fixed costs associated with hiring regular drivers.}

Overall, delivery options provide more flexibility but, at the same time, they increase the planning complexity, both when they refer to location options (multiple locations are offered) or transportation options (multiple means of transportation or different cohorts of drivers are available). Next, we first focus on delivery location options in Section \ref{sec:where}, including two MILP models with deterministic input data in Section \ref{formu_option} and the related literature in Section \ref{lite:where}. Then,  we discuss the literature related to transportation options in Section \ref{sec:deliveryMean}.

\subsection{Where to deliver}
\label{sec:where}
In this section, we focus on the problems where different options are available as to where the parcels could be delivered.
\subsubsection{Mathematical formulation}
\label{formu_option}

We present a MILP model based on the setting of the problem studied in \cite{ulmer2019same} with some modifications. Specifically, as in the former sections, we consider a single uncapacitated vehicle and we assume to have perfect information, while in \cite{ulmer2019same} orders are unknown. Also, the time at which orders are placed is known beforehand. 
The problem is described as follows.

A set of customers $N$ orders goods online over a discrete time horizon $T$ = $(0, 1, ...,  T_E)$. We denote the order time {(release date)} of each customer with $r_n$, $n \in N$. Deliveries are not made to customers' locations but to a set of pickup stations, where multiple parcels can be delivered in a single trip. Specifically, we denote as $\mathcal{P}$ the set of pickup stations and as 0 the depot. $\mathcal{P}_0$ is the union of $\mathcal{P}$ and $0$.
{We assume that each customer 
can be served from any 
of the stations that are within a travel radius centered around the customer's location.} We denote as $P_n$ the set of stations that can be used to serve customer $n$. The delivery is performed by a single uncapacitated vehicle. Consistently with the assumption made in \cite{ulmer2019same}, we assume that the vehicle only {travels from the depot to a pickup station and back. This type of trip is called direct trip}. Thus, a trip corresponds to a visit to a pickup station, and no routing decision is involved.   The capacity of pickup station $i \in \mathcal{P}$ during the entire time horizon is $c_i$, and the newly 
freed space due to orders being picked up is not considered. Each parcel occupies one unit of station capacity. The travel time from the depot to a pickup station $i$ and from a pickup station $i$ to a depot is given by $d_{0i}$ and $d_{i0}$ respectively, $i \in \mathcal{P}$.  A solution consists of up to $|\cK|$ trips performed by the vehicle, where the number of $|\cK|$  {is not greater than $|N|$}. The objective is to minimize the sum of delivery times, i.e., the sum over all customers of the difference between the \emph{earliest pickup time} (the time when the order is delivered to a pickup station) and the order time. In order to serve as many customers as possible before the deadline $T_E$, a big penalty value is added to the earliest pickup time in case a customer is not served.

The mathematical formulation of the problem makes use of the following decision variables:
\begin{enumerate*}[(i)]
    \item $y_{ni}^k$: binary variable equal to 1 if the order of customer $n$ is served in trip $k$ to pickup station $i$, 0 otherwise;
    \item $y_n$: binary variable equal to 1 if customer $n$ is served, 0 otherwise;
    \item $s_k$: continuous variable representing starting time of trip $k$;
    \item $t_n$: continuous variable representing earliest pickup time of order $n$;
    \item $z_i^k$: binary variable equal to 1 if pickup station $i$ is visited in trip $k$, 0 otherwise. 
\end{enumerate*}
The formulation is the following: 
\small
\begin{subequations}
	\label{eq:base3}
\begin{equation}\label{eq:obj_ILP3_base}
\hbox{min} \sum_{n \in N} (t_n - r_n)
\end{equation}
\begin{align}
\label{eq:each_ILP3_most_once1}
 \sum_{i \in P_n}\sum_{k \in \cK}y^k_{ni}  &= y_n && \quad n \in N \\
\label{eq:y0_ILP3_base}
 \sum_{i \in \mathcal{P}\backslash P_n}y^k_{ni} &= 0 && \quad n \in N, k \in \cK\\
\label{eq:loca_ILP3_capacity}
    \sum_{k \in \cK}\sum_{n \in N}y_{ni}^k & \leq c_i && \quad i \in \mathcal{P}\\
\label{eq:Y_ILP3_base}
     t_n &\geq M(1-y_n) && \quad n \in N\\
\label{eq:loca_ILP3_trip1}
    c_iz_i^k & \geq \sum_{n \in N}y_{ni}^k && \quad i \in \mathcal{P}, k \in \cK\\
\label{eq:non-necessaryz_ILP3}
    \sum_{n \in N}y_{ni}^k & \geq z_i^k && \quad i \in \mathcal{P}, k \in \cK\\
\label{eq:start_ILP3_base}
    s_k & \geq \sum_{i \in P_n} r_n 
    y_{ni}^k   && \quad n \in N, k \in \cK\\
\label{eq:timeHorizon_ILP3_base}
    s_{|\cK|} & \leq T_E && \\
\label{eq:next_ILP3_base}
    s_{k+1} & \geq s_k + {\sum_{i \in \mathcal{P}}(d_{0i}+d_{i0})z_i^k} && \quad k \in \cK\\
\label{eq:p_ILP3_base}
     t_n & \geq s_k + \sum_{i \in P_n}d_{0i}y_{ni}^k - (1-\sum_{i \in P_n}y_{ni}^k)M && \quad n \in N, k \in \cK\\
\label{eq:loca_ILP3_trip2}
    \sum_{i \in \mathcal{P}}z_i^k & \leq 1 && \quad k \in \cK\\
    y^k_{ni}, z_i^k, y_n & \in \{0,1\} && \quad n \in N, i \in \mathcal{P}_0, k \in \cK
    \\
    s_k, t_n &\geq 0 && \quad n \in N, k \in \cK.
\end{align}
\end{subequations}
\normalsize
The objective function \eqref{eq:obj_ILP3_base} minimizes the service time over all customers. Constraints (\ref{eq:each_ILP3_most_once1}) limit each customer to be served at most once, and  (\ref{eq:y0_ILP3_base}) ensure that customers are not assigned to pickup stations outside the radius. Inequalities  (\ref{eq:loca_ILP3_capacity}) establish the maximum capacity for each pickup station while (\ref{eq:Y_ILP3_base}) associate a large penalty value with customers who are not served.  Constraints (\ref{eq:loca_ILP3_trip1}) state that a customer can be served by a station in a trip only if there is a delivery to the station and the capacity is not exceeded. Inequalities (\ref{eq:non-necessaryz_ILP3}) guarantee that a pickup station is not visited in a trip when no parcel is delivered to the station and (\ref{eq:start_ILP3_base}) guarantee that each trip can only serve customers after order time. Constraints (\ref{eq:timeHorizon_ILP3_base}) state that the last trip must finish before the deadline and (\ref{eq:next_ILP3_base}) set the starting time of each trip with respect to the previous trip.  Constraints (\ref{eq:p_ILP3_base}) set the earliest pickup time for each customer where $M$ is a big value. Finally, (\ref{eq:loca_ILP3_trip2}) impose that a single pickup station can be visited in each trip.

The difficulty of the problem lies in satisfying simultaneously the capacity of pickup stations and the deadline. In a more general case, instead of considering direct trips only, the formulation can be adapted to include routing among multiple pickup stations in the same trip. The new problem shares similarities with the Generalized Traveling Salesman Problem (see \cite{laporte1996some}). The formulation is presented below.

First, we introduce an additional class of decision variables $a_i^k$ representing the arrival time at the pickup station $i$ in trip $k$,  $i \in \mathcal{P}$. In addition,  $x_{ij}^k$  are binary variables equal to 1 if arc $(i, j)$ is traversed in trip $k$, with $i, j \in \mathcal{P}_0$. We denote as $A$ the set of arcs formed by any two nodes $i, j$ in $\mathcal{P}_0$.  {Note that, in our model, variable $z_i^k$ is defined not only for pickup stations but also for the depot, i.e., $i \in \mathcal{P}_0$.} The new formulation reads as follows. 


\small
\begin{subequations}
	\label{eq:base4}

\begin{equation}\label{eq:obj_ILP4_base}
\hbox{min} \sum_{n \in N} (t_n - r_n)
\end{equation}
 \begin{align}
\text{  $\qquad (Y_n,y^k_{ni},z^k_i,s_k,t_n)$ satisfies \eqref{eq:each_ILP3_most_once1}-\eqref{eq:timeHorizon_ILP3_base}
} \notag 
\end{align}
\begin{align}
\label{eq:inout_ILP4_base}
 \sum_{(i, j)\in A}x^k_{ij} &= \sum_{(j, i)\in  A}x^k_{ji}  = z_i^k && \quad i \in \mathcal{P}_0, k \in \cK\\
\label{eq:next_ILP4_base}
    s_{k+1} & \geq s_k + \sum_{(i, j) \in A}d_{ij}x_{ij}^k && \quad k \in \cK\\
\label{eq:sa_ILP4_base}
    s_k & = a_0^k && \quad k \in \cK\\
\label{eq:arrivaltime1_ILP4_base}
    a_j^k & \geq a_i^k + d_{ij}x_{ij}^k - (1-x_{ij}^k)M && \quad i\in \mathcal{P}_0, j \in \mathcal{P}, k \in \cK\\
\label{eq:arrivaltime2_ILP4_base}
     t_n & \geq a_i^k-(1-y_{ni}^k)M && \quad n\in N, i \in P_n, k \in \cK\\
    x^k_{ij}, y^k_{ni}, z_i^k, y_n & \in \{0,1\} &&\quad i,j \in \mathcal{P}_0, k \in \cK \\
    s_k, t_n,  a^k_i &\geq 0 &&\quad i,j \in \mathcal{P}_0, k \in \cK 
\end{align}
\normalsize


Constraints (\ref{eq:inout_ILP4_base}) are degree constraints and state that a parcel can be delivered to a station only when the station is visited.  Constraints (\ref{eq:next_ILP4_base}) fix the starting time of consecutive trips while (\ref{eq:sa_ILP4_base}) set the starting time of trip $k$ as the time at which the vehicle leaves from the depot. Inequalities (\ref{eq:arrivaltime1_ILP4_base}) link the arrival times of nodes visited consecutively in a trip and (\ref{eq:arrivaltime2_ILP4_base}) link pickup times with arrival times. Similarly as in Formulation \eqref{eq:base2}, subtours are eliminated through (\ref{eq:arrivaltime1_ILP4_base}) and subtour elimination constraints (\ref{eq:subtour_ILP4_base}) 
\small
\begin{align}
    \sum_{i \in S, j\in \mathcal{P}_0\setminus{S}}x_{ij}^k \geq z_l^k, && S \subseteq  \mathcal{P}_0, 0 \in S, l \in \mathcal{P}\setminus{S}, k \in \cK&&\label{eq:subtour_ILP4_base}
\end{align}\end{subequations}  \normalsize
can be introduced to strengthen the formulation.

\bigskip

As for the former cases, formulations (\ref{eq:base3}) and (\ref{eq:base4}) can be adapted to solve similar problems. For instance, instead of minimizing the service time, the objective function can be maximizing the number of requests served before the deadline {(which is anyway implicitly considered in the current formulations with fixing the value of $M$). Also, one might want to consider the case without a deadline and where all customers have to be served, by searching for the minimum cost solution.}



\subsubsection{Literature review on delivery location options in SDD}
\label{lite:where}
Providing various delivery options is a way to reduce delivery failures. The related research problems span from determining the locations and layouts of pickup stations (\cite{kahr2022determining}) to selecting a delivery place for each request 
\citep{zhou2018multi, dumez2021large, tilk2021last, mancini2021vehicle, dragomir2022pickup}.  {There are different contributions to problems involving multiple delivery locations. However, the ones considering dynamic customer arrivals in SDD services are still limited.} 

\cite{ulmer2019same} is the first study that combines pickup stations and autonomous vehicles for SDD. Each order is associated with a preferred pickup station, and adjacent stations (called neighborhood) within certain customer travel time around the preferred station are also accepted with no penalty.
The problem is to decide where and when to dispatch a vehicle and which parcels to deliver on each dispatch. Autonomous vehicles perform deliveries from the depot to the stations. The impact of customer behavior, speed of autonomous vehicles, and capacities of vehicles and pickup stations are considered. The study addresses the tradeoff between consolidation and fast dispatches by defining a threshold, which specifies the minimum number of parcels that can be shipped in each dispatch. To determine the threshold, the authors propose a set of policy function approximation (PFA) approaches.
A sample average approximation over 1000 scenario realizations is used to find the policy that minimizes service time. {In the computational section, the authors analyze three types of neighborhoods: none, normal, and large, where none means each order is delivered to the preferred pickup station and large means a larger neighborhood size than a normal one is eligible for delivery. They find that by increasing neighborhood size, the percentage of deliveries to the preferred stations decreases but the service level improves.}

\subsection{Literature review on how to deliver}
\label{sec:deliveryMean}
Drones, a recent technology that is starting to be used in deliveries, offer several benefits with respect to vehicles, the main ones being that they are autonomous and they are not impacted by  (and do not impact) 
the {road} traffic. {Drones can fly for a limited time, but they can still travel back and forth in middle-sized cities within 15km\footnote{https://www.cnet.com/tech/computing/amazon-drone-deliveries-coming-to-texas-city-as-prime-air-expands/}. When we compare vehicles with drones, we notice that vehicles are capable of delivering many orders at once but often suffer from traffic congestion, while the capacity of drones is limited but they can travel at higher speed (as they avoid traffic) and with free trajectories.} However,  drones are sensitive to other factors, such as the wind conditions in landing zones and the density of the population \citep{gianfelice2022real, weibel2005safety}. 
Likewise, following the concepts introduced in Section \ref{classification_scheme},  we summarize the papers presented in Sections \ref{lite:where} and \ref{sec:deliveryMean} to give an overview and provide features of Models (\ref{eq:base3}) and (\ref{eq:base4}) as a point of reference in Table \ref{table:literature_deliveryOption}. We see that there is high interest in studying the possibility of using drones to resupply vehicles or transshipment points for transferring goods.

\begin{table}[!h]
\caption{Summary of papers on delivery options. The abbreviations: heterogeneous fleet (HF), multiple trips per vehicle (MT), maximum route duration (RD), order pickup range (PR), vehicle capacity (VC), time windows (TW), delivery deadlines (DD), battery charging (BC), and pickup station capacity (PS).}
\label{table:literature_deliveryOption}
\scalebox{0.62}{\color{black}
\begin{tabular}{p{0.2\textwidth}p{0.17\textwidth}p{0.2\textwidth}p{0.2\textwidth}p{0.16\textwidth}p{0.18\textwidth}p{0.2\textwidth}p{0.13\textwidth}}\hline
\textbf{Literature}     & \textbf{Problem}    & \textbf{Objective}                                                                          & \textbf{Uncertainty}                                                                       & \textbf{Uncertainty Nature} & \textbf{Decision}              & \textbf{Constraint} & \textbf{Solution Approach}            \\ \hline
\cite{ulmer2019same}       & vehicle dispatching                       & min(sum of delivery times)                                                                                                                  & requests(arrival time), customer pickup   time                                                & dynamic \& stochastic & schedule; delivery locations                        & DD; MT; PS; PR; VC  & ADP                            \\
\cite{ulmer2018same}          & VRP                                       & max(number of requests served)                                                                                                              & requests(arrival time)                                                              & dynamic \& stochastic & reject/accept; assignment of services; route & DD; MT; HF; BC      & ADP                            \\
\cite{chen2022deep}              & VRP                                       & max(number of requests served)                                                                                                              & requests(arrival time, location)                                                              & dynamic \& stochastic & reject/accept; assignment of services; route & DD; MT; HF; BC      & RL                             \\
\cite{dayarian2020same}          & VRP with drone resupply                 & max(percentage of orders served); min(operational costs)                                                                                    & requests(arrival time)                                                                        & dynamic \& stochastic & schedule; route; drone resupply                  & TW; DD; MT; HF; VC & two-phase heuristic            \\
\cite{mccunney2019simulation} & vehicle dispatching with drone resupply & min(travel time)                                                                                                                            & requests(arrival time)                                                                        & dynamic \& deterministic                & schedule; locations of transshipment points               & DD; MT; HF; VC      & MILP; simulation               \\
\cite{pina2021traveling}         & TSP (with drone resupply)                 & min(completion time for serving all customers)                                                                                              &                              no                                                             &      static \& deterministic                  & route; drone resupply                  & HF; MT; VC          & two-stage decomposition heuristic                 \\
\cite{lin2018demand}                & VRP                                       & min(sum of transportation cost, fuel cost, and emissions cost)                                                                      & requests(arrival time, pickup location, delivery location)                                    & dynamic \& deterministic               & route                            & DD; VC                  & MILP; heuristic \\ 
\cite{arslan2019crowdsourced}             & matching and TSP            & min(distribution cost: sum of compensation paid to the drivers and cost of dedicated vehicle trips) & requests(arrival time), driver (arrival time)                                                 & dynamic \& deterministic                & match; route                     & TW; VC; RD; HF          & rolling horizon       \\
\cite{archetti2021online}           & VRP                                       & min(distribution cost)                                                                                                                      & requests(arrival time)                                                                        & dynamic \& deterministic                & route                            & TW; VC; DD; HF          & local search        \\ 
\cite{ausseil2022supplier}            & matching   & max(expected revenue)   & request (arrival time, pickup location, delivery location), driver (arrival time, choice) & dynamic \& stochastic & match                            & DD; VC                  & scenario sampling  \\\hline
Model (\ref{eq:base3})  & OP & min(service time) & no                                                 & static \& deterministic         & schedule; delivery locations;                 & DD; SV; MT; PS         & MILP   \\
Model (\ref{eq:base4})  & OP & min(service time) & no                                                 & static \& deterministic         & schedule; delivery locations; route                & DD; SV; MT; PS          & MILP                               \\\hline
\end{tabular}}
\end{table}

\cite{ulmer2018same} present a dynamic vehicle routing problem with a heterogeneous fleet of drones and vehicles.  Customers order goods during the day and they have to be delivered by the end of the day. 
{The authors propose to use a threshold for dividing the service
area into two zones. Then, following the observation that vehicles may be more suitable to serve customers located in downtown areas close to the depot while drones may be more appropriate for rural areas with dispersed customers, vehicles are preferred to serve customers in the zone within the threshold.
Drones preferably serve the ones outside. The authors then propose a PFA parameterized over such a threshold. Among a set of potential threshold values, they select the PFA leading to the highest average number of customers served within the delivery deadline.} Compared to other threshold-type benchmark policies, \cite{ulmer2018same}  find that districting by the proposed approach is highly profitable as it considerably increases the expected number of requests served.

\cite{chen2022deep} also make the choice of delivery means between vehicles and drones. Rather than deciding on the basis of the locations of customers only (as done in \cite{ulmer2018same}), the authors use more information such as time and fleet availability to decide whether to use a vehicle or a drone. They propose a deep Q-learning approach to determine the assignment of a new customer to either a vehicle or a drone. The computational results compared to benchmark policies show the effectiveness of the proposed approach. 

There is a vast literature on combining drone deliveries with conventional trucks and the survey on this topic would require an article on its own. We refer interested readers to recent survey articles by \cite{chung2020optimization}, \cite{macrina2020drone}, and \cite{moshref2021applications}. When it comes to truck and drone routing dealing with online orders in SDD, the number of contributions is rather limited. {\cite{dayarian2020same} propose a home delivery system in which
a fleet of delivery trucks is dispatched from a fulfillment center to perform home deliveries of online orders. The trucks are ``collaborating'' with a fleet of drones 
which is used to regularly resupply them with new parcels. That way, 
the vehicles do not need to return to the depot for serving new customers. 
In the beginning, the trucks are dispatched to serve a predetermined route. Drones can be dispatched to deliver additional parcels to a truck at a pre-determined meeting location and time.
In that case, the route of the truck is re-optimized to guarantee that the parcels which are on board are delivered on time. 
In order to decide how the truck should be resupplied, the authors propose two strategies: restricted resupply and flexible resupply. 
For the former, resupply can happen only once all on-board orders are
delivered, whereas in the latter case, the truck may be instructed to perform a detour towards a meeting location where it will be resupplied with new parcels. The computational results show that the flexible resupply performs slightly better than the restricted resupply.}

\cite{mccunney2019simulation} also study the impact of using drones to resupply trucks to avoid trucks returning to a depot for picking up new orders. They propose a MILP model to decide the meeting locations of drones and trucks and the customers assigned to each of these locations. By doing simulation tests based on location and delivery data from two American cities, they find that using drones to resupply trucks can reduce delivery times and distances for serving all customers compared to a traditional truck-only delivery system.

\cite{pina2021traveling} study a similar problem using drones to resupply trucks as done in \cite{dayarian2020same} and \cite{mccunney2019simulation}, but they assume release dates and locations of orders are known in advance. They model the problem as the traveling salesman problem with release dates and drone resupply. First, they provide a MILP formulation with the goal of minimizing the time to deliver all orders. Then they propose a two-stage decomposition heuristic: the first stage determines a truck route; the second stage defines the locations where each order will be loaded onto the truck by the drone. The experiments compare the system to a traditional truck-only delivery system and show that using drone resupply can significantly reduce the total delivery time.

The concept of crowdshipping is based on {incentivizing
private citizens to deliver packages by offering them a compensation}. This strategy is widely 
deployed nowadays (see  Walmart (\cite{WalMart} and Amazon (\cite{Bensinger}) for early implementations of the idea). \cite{ulmer2020workforce} and \cite{dai2020workforce} investigate workforce scheduling and capacity planning of crowdsourced drivers respectively to achieve a satisfying service level at minimum cost.
\cite{dayarian2020crowdshipping} quantify the possible benefits of crowdshipping for SDD and demonstrate the value of future information in a highly dynamic and stochastic SDD environment. 

When it comes to comparing different delivery paradigms, \cite{lin2018demand} consider three options: hub-and-spoke, SDD with a commercial fleet, and SDD with crowdshipping.  
In the first one, all requests are collectively received by a hub-and-spoke delivery system by the end of the day and scheduled for next-day service. In the two SDD paradigms, requests are added in real-time and are to be served during the same working day. The results find hub-and-spoke to be cost-effective for the traditional delivery service but expensive for the SDD service. They also find that crowdshipping provides a low-cost solution to SDD.

Routing problems with crowdshipping are also called \textit{routing problems with occasional drivers}. This term was introduced in \cite{archetti2016vehicle}. In the context of SDD, we highlight several contributions. 
\cite{arslan2019crowdsourced} investigate a problem in which a  provider receives orders online and matches them to ad hoc drivers or dedicated vehicles 
for which the associated routes are also determined. 
Orders and drivers  arrive dynamically over time. 
The  objective is to minimize the sum of the compensation paid to the drivers matched and the cost of the dedicated vehicle trips to serve all customers. To solve the matching problem, the authors propose a rolling horizon framework and develop an exact solution method that is applied whenever there is new information. 

\citet{archetti2021online} also consider a setting in which customer orders arrive online. The company managing deliveries has a set of ``regular'' drivers (drivers working regularly for the company) as well as crowdsourced drivers (occasional drivers). {There are four major differences compared to the study of \cite{arslan2019crowdsourced}:  the information related to occasional drivers is deterministic, not all parcels need to be delivered on the same day, all parcels must be dispatched from the same depot, and preemptive depot returns are allowed.} The objective is to serve all customers at a minimum cost. Besides the routing cost for regular drivers and the compensation cost for occasional drivers, the authors also introduce penalties for violating time windows or postponing the delivery to the next day. An iterated local search is proposed to solve the problem which outperforms an insertion heuristic.

Finally, \cite{ausseil2022supplier}  consider a peer-to-peer logistics platform to match crowdsourced occasional drivers dynamically. Contrary to what is assumed in the contributions mentioned above, crowdsourced drivers might not accept serving a customer. In order to mitigate acceptance uncertainty, the platform offers each driver a menu of requests to choose from.  {The objective is to maximize the expected revenue, 
{i.e., the platform rewards obtained from the successful assignments of orders to drivers} within the time horizon.} To treat the menu selections, the authors present a multiple-scenario approach, including repeatedly sampling potential driver selections, solving two-stage stochastic models, and deriving solutions through a consensus algorithm.

\section{Pickup and delivery}
\label{pickup_delivery}
When parcels are associated with origins and destinations, the problems are known as pickup and delivery problems (\cite{savelsbergh1995general}). 
In local delivery networks, pickup and delivery operations are typically combined so that drivers can be assigned to both operations. The typical objective remains the minimization of the routing cost. The difference with respect to the problems analyzed in the former sections is that parcels have multiple origins (instead of a single one). In same-day services, however, not all customers are known in advance, and the uncertainty can lead drivers to revisit previously served neighborhoods to deliver new parcels. A family of same-day pickup and delivery problems has emerged recently, and this section highlights some of them. We introduce the literature on same-day pickup and delivery problems in Section \ref{lite:pure_pick}, meal delivery service in Section \ref{meal_delivery}, and personal shopper service in Section \ref{lite:personal}. The papers are summarized in Table \ref{table:literature_pickup}. We observe that performance measures are more varied than those in previous sections. Additionally, there are new uncertainty factors such as driver parking time and food ready time. Furthermore, there is a greater diversity in solution approaches.

\begin{table}[!h]
\caption{Summary of papers on pickup and delivery. The abbreviations: machine learning (ML), reinforcement learning (RL), multiple trips per vehicle (MT),  vehicle capacity (VC), time windows (TW), and delivery deadlines (DD).}
\label{table:literature_pickup}
\scalebox{0.63}{\color{black}
\begin{tabular}{p{0.2\textwidth}lp{0.23\textwidth}p{0.2\textwidth}p{0.16\textwidth}p{0.24\textwidth}lp{0.13\textwidth}}\hline
\textbf{Literature}     & \textbf{Problem}    & \textbf{Objective}                                                                          & \textbf{Uncertainty}                                                                       & \textbf{Uncertainty Nature} & \textbf{Decision}              & \textbf{Constraint} & \textbf{Solution Approach} \\ \hline     \cite{mitrovic2004waiting} & VRP      & min(travel distance)                     & requests(arrival time, pickup \& delivery   locations)                                                                    & dynamic \& deterministic                & schedule; route                                                                                & TW; DD     & cheapest insertion; waiting strategies      \\
\cite{pureza2008waiting}           & VRP      & min(\#lost sales, routes, travel distance)                   & travel time, requests(arrival time)                                              & dynamic \& deterministic                & accept/reject; route; schedule                                                                & TW; DD; VC & waiting strategy; request buffering strategy; insertion heuristics  \\
\cite{ghiani2009anticipatory}                & VRP      & min(customer inconvenience)                                  & requests(arrival time)                                                                                                 & dynamic \& stochastic & schedule; route               & TW; DD;    & anticipatory insertion; local search  \\
\cite{ackvaconsistent}             & VRP      & max(\#requests served)                                       & requests(release time, pickup store, location) & dynamic \& stochastic                 & route; schedule                                                                                & TW; DD; VC & scenario sampling \\
\cite{reyes2018meal}                 & VRP      & optimize performance measures (\#requests served, courier compensation, cost per order...)                               & requests(arrival time)& dynamic \& deterministic                & assignment of services; route; schedule & DD; MT     & rolling horizon\\
\cite{yildiz2019provably}       & VRP      & optimize performance measures (courier compensation, click-to-door time...)                                &                                                 no                                      & static \& deterministic                 & assignment of services; route; schedule                                          & DD; MT     & column and row generation    \\
\cite{ulmer2021restaurant}                 & VRP      & min(\#delays)                                                & requests(arrival time), food ready time                                                                                & dynamic \& stochastic & route; assignment of services                                                                             & TW; DD; MT & anticipatory customer assignment policy                                                                                                                                                                                                                      \\
\cite{mao2019faster}                   & VRP      & no   &     no &     static \& deterministic &     assignment of services  &     TW      & empirical study     \\
\cite{hildebrandt2021supervised}       & VRP      & max(\#requests served) & requests(restaurant choice), food ready time, driver (parking and   delivery times) & dynamic \& stochastic & route                                                                                         & TW; DD     & ML                                \\
\cite{auad2021courier}                  & matching & balance courier and customer satisfaction                                   & requests(arrival time, ready time), delivery (capacity)                                                                & dynamic \& deterministic                & match                                                                                          & TW; DD     & rolling horizon    \\
\cite{arslan2021operational}                & VRP      & max(\#requests served)                                       & requests(arrival time)                                                                                                 & dynamic \& deterministic               & accept/reject; route                                                                          & TW; DD; VC & three order-consolidation strategies                           \\ \hline                            
\end{tabular}}
\end{table}

\subsection{Literature review on same-day pickup and delivery problems}
\label{lite:pure_pick}

One of the first contributions studying a {same-day pickup and delivery problem (SDPD)} is due to \cite{mitrovic2004waiting}. The authors focus on same-day pickup and delivery of letters and small parcels by investigating the impact of deploying four possible waiting delivery strategies. The first strategy, called the drive-first waiting strategy (DF), requires a vehicle to leave its current location at its earliest departure time. The second, named wait-first strategy (WF), on the opposite, requires a vehicle to leave at its latest possible departure time. The next two are combinations of DF and WF. In the third, named dynamic waiting strategy (DW), each vehicle drives as soon as possible while serving close locations; when all such locations are visited and the vehicle has to serve the first distant location, it waits as long as possible. The last strategy, called advanced dynamic waiting strategy (ADW), propagates the waiting time available in the route along the whole route, i.e., assigning it to different locations visited in the route.  Based on the tests with real-life data, the authors show that WF builds shorter routes than DF but uses a larger number of vehicles, and ADW is the most efficient with respect to total route length and the number of vehicles used. 

\cite{pureza2008waiting} study a similar problem as \cite{mitrovic2004waiting} but they allow  rejection of requests. Also, the distances are time-dependent, i.e., they may vary according to the time of the day. The objective is to determine a set of routes minimizing hierarchically the number of lost sales (requests rejected upon their arrivals), the number of routes, and the total travel distance within the time horizon. They propose two delivery strategies. {The first one is a waiting strategy, which keeps the vehicle waiting at the current location for some time to accumulate requests with new arrivals rather than serving the next planned location immediately after the service of the current location. The second is a request buffering strategy, which schedules route adjustments to include new but non-urgent requests at a later stage instead of adjusting the route upon each new arrival.} {The numerical results show that using the waiting strategy is more advantageous compared to the request buffering in terms of lost requests and the number of vehicles required. However, the request buffering strategy results in a shorter route distance compared to the waiting one.} 

Another early contribution is due to \cite{ghiani2009anticipatory}. The authors investigate a multi-vehicle dispatching problem with pickup and delivery where the requests arrive according to a known stochastic process. 
A sampling technique, which is used to anticipate future demands, is embedded in the cheapest insertion heuristic and a local search algorithm to design the delivery routes. 

The studies of \cite{lin2018demand} and \cite{arslan2019crowdsourced} (introduced in Section \ref{sec:deliveryMean}) also address SDPD.
When comparing three delivery paradigms, \cite{lin2018demand} also consider delivery requests that are associated with a pair of pickup and delivery locations. 
Similarly, when matching the orders with ad-hoc drivers and dedicated vehicles, \cite{arslan2019crowdsourced} allow orders to be associated with different pickup locations (a warehouse, a dedicated pickup point, or a retail shop). 

\cite{ackvaconsistent} {introduce} a pickup and delivery problem  with micro-hubs.
The pickup and delivery locations of each parcel are mapped to their nearest micro-hubs.
Hence, different parcels may share the pickup and delivery locations, and vehicles can exploit consolidation opportunities by bundling parcels from different pickup locations for joint delivery to the same area. Vehicles follow a predetermined daily schedule visiting micro-hubs at predetermined times and in a predetermined sequence. 
The objective is to find a consistent routing schedule, i.e., visiting the same micro-hubs at the same times every day, that maximizes the expected number of parcels delivered daily.
The authors model the problem as a two-stage stochastic program: the first stage determines a sequence of stops at micro-hubs with departure times; the second stage decides which parcels to serve given the schedule. The problem is solved using scenario sampling as introduced by \cite{bent2004scenario} and four problem-specific consensus functions (CF) are proposed. 
The computational results show their approach outperforms the original CF proposed by \cite{bent2004scenario} and two benchmark heuristics.

\subsection{Literature review on meal delivery service}
\label{meal_delivery}

{Restaurant meal delivery is a typical application for same-day/instant pickup and delivery problems.  Fierce competition in this sector requires efficient scheduling and routing optimization algorithms, in order to gain a competitive advantage. In these applications, we distinguish between the \emph{order time} and \emph{ready time}. 
Order time corresponds to the time when a customer sends a request, choosing from various restaurants. Ready time is the time when the ordered meal becomes available at the restaurant.
The customers are usually unknown until orders are placed, and the ready time of food is unknown as well, which creates tremendous challenges in managing the delivery service. Once a request is sent, the order is passed to a delivery company, which matches it to a courier for picking up and delivering. The delivery company has to anticipate future uncertain requests, in terms of location and timing, as well as the ready time and delivery routes, before deciding how to match orders with the couriers.  }

{\cite{reyes2018meal} is the first article to introduce a model that incorporates the main structural features of a meal delivery routing problem, such as:  multiple restaurants as pickup points, orders arriving dynamically, delivery capacity distributed around the day in the way of courier shifts, or order bundling opportunities. The authors assume the information about restaurants and couriers is known in advance but each order is known only at its order time. The problem is to determine feasible routes for couriers to complete the pickup and delivery of orders. The objective is to optimize single or multiple performance measures, such as click-to-door time (i.e., the delivery time), the number of orders delivered, or the cost per order (i.e., total courier compensation divided by the number of orders delivered). 
To deal with the uncertainty, the authors propose a rolling horizon matching-based algorithm that 
matches one or a bundle of orders with a courier. After the matching of orders, a postponement of less time-critical orders is applied. The results indicate that the proposed method can solve the set of benchmark instances in (near) real-time and that the obtained solutions are of high quality with respect to the performance measures mentioned above. } 


\cite{yildiz2019provably} develop an exact solution approach for the meal delivery problem introduced by \cite{reyes2018meal}. The authors assume perfect information concerning the {order arrival stream including order request time, order-placement time, order-ready time, and the drop-off location}.
The solution approach is based on simultaneous column and row generation. An extensive numerical study using instances from real-world historical data is conducted. {The authors derive several interesting managerial insights, showing that cost and service objectives are well aligned and that the minimum payment to couriers can be guaranteed without a significant increase in cost. }

{\cite{ulmer2021restaurant} study a similar problem setting as \cite{reyes2018meal}, 
but they assume that the probability distributions of {order times,} {ready times}, and customer locations are known. A delivery company assigns each order to a courier who picks up and delivers the order. Multiple orders can be assigned to the same courier, and the major challenge is  in having accurate estimations of 
ready times, in order to prevent waiting of the courier and, thus, potential delays along the route. The customers are promised that the orders are delivered by a certain deadline and the overall objective is to minimize the expected delays.
The authors model the problem as an MDP and propose an anticipatory customer assignment policy with a time buffer as a type of cost function approximation (\cite{powell2019unified}). {The time buffer is a ``cost'' added to the scheduled arrival times at the customers to account for the stochasticity.} 
The policy postpones the assignment of  a courier to an order to get more information on ready time and future orders. As in \cite{reyes2018meal}, the postponement enables the potential bundling, i.e., assigning  multiple orders at once to the same courier. In the computational study, \cite{ulmer2021restaurant} demonstrate the importance of the proposed time buffer compared to solely using postponement and bundling as in \cite{reyes2018meal}.} 

{While the studies mentioned above assume fixed customer choices, they do not consider the possible impact of late deliveries on customer satisfaction.} \cite{mao2019faster} examine the importance of on-time and early delivery {empirically based on two-month transactional data from an online meal delivery platform in China}. They find that the negative impact on customer satisfaction due to late deliveries is much more severe than the positive impact created by early deliveries.

\cite{hildebrandt2021supervised} also work on meal pickup and delivery problems. A novelty of the paper lies in considering that customer choice is affected by estimated arrival time. That is to say, customers can choose between different products and services, and their choice is based on the restaurant preference and the expected time of delivery. 
Thus, the authors work on estimating meal arrival times to inform customers at the time they make their selection. Two supervised learning approaches for the prediction of meal arrival times are proposed. 
The results show that accurate arrival times can raise service perception and guide customers' selection to improve the delivery system. Besides, the authors find that arrival time predictions help reduce delays but may block customers from ordering.

Besides satisfying customers, another aspect to be considered is courier satisfaction. {The works mentioned above all assume that couriers do not execute any autonomous decisions.} In fact, couriers may reject certain deliveries due to long distances. Logistics companies in the field typically compete for the same set of couriers who can subscribe to multiple platforms. Thus it is necessary to maintain courier satisfaction as well. \cite{auad2021courier} consider a pickup and delivery problem by combining both customer-centric and courier-centric performance metrics. Assuming couriers prefer to operate in relatively small geographic areas, they propose an optimization-based rolling horizon algorithm for managing couriers by coping with region resizing and delivery task assignment decisions. {Their empirical study using real-world data shows that for a fixed courier fleet size, courier satisfaction metrics can be increased without sacrificing key performance metrics related to customer service.}

\subsection{Literature review on personal shopper service}
\label{lite:personal}
Personal shopper service experienced significant growth in online grocery (\cite{personalShopper}). Compared to delivering fresh meals, there is more time flexibility here. However, the personal shopper service often requires visiting multiple stores to serve a single customer {rather than visiting one pickup location per customer. Thus, the problems can be seen as an extension of pickup and delivery problems with multiple pickup locations per order}. 

Pickup and delivery problems in on-demand personal shopper service are investigated by \cite{arslan2021operational}.   
In a classical pickup and delivery setting, 
the pickup time is often negligible compared to travel times. However, in personal shopper service applications, the pickup time cannot be ignored, as it increases with the number of simultaneous 
pickups per store visit. 
\cite{arslan2021operational} aim to maximize the number of requests served, subject to fixed shopper capacity (which is the number of tasks a shopper can do). The authors propose three different order-consolidation strategies for service operators. The obtained results show that order consolidation, especially shopping consolidation, can significantly increase the number of requests served within the available shopper capacity.

\section{Districting - fleet size and mix}
\label{Districting}


{Decisions at distribution companies are carried out at three levels (see, e.g., \cite{jabali2012continuous}).  At the \emph{strategic} level, the decisions are made concerning the types and numbers of vehicles that need to be purchased to fulfill the expected demand. 
At the \emph{tactical} level, leasing additional capacity or subcontracting part of the fleet might be needed. Finally, at the \emph{operational} level, daily operations, which usually concern the routing of vehicles, are conducted. 
So far, we focused on operational problems, however, there are other relevant problems at the strategic and tactical levels that need to be considered in the context of SDD. 
One of the first optimization problems that integrate strategic and tactical decisions in the general context of routing is introduced by \cite{golden1984fleet}. Their \emph{fleet size and mix vehicle routing problem} asks to determine the optimal number of heterogeneous vehicles with possibly different costs and capacity to lease or purchase. The authors consider the problem of minimizing the cost including vehicle-dependent fixed cost and variable cost components such as fuel, maintenance, and manpower.} 

{In SDD, a single depot tends to serve a 
relatively large region with hundreds or even thousands of parcels delivered every day. Hence,  it is highly impractical to 
assume each vehicle should be able to serve the whole region.
Instead, customer locations are partitioned
into zones (districts) with dedicated vehicles, where
each district can be served by either a single vehicle (\cite{banerjee2022fleet}) or a small fleet (\cite{stroh2022tactical}). 
Such partition of a service region is also known as  \emph{districting}.   
Thus, the districting problems belong to the strategic or tactical decision-making levels (\cite{haugland2007designing}). Next, we will first discuss in Section \ref{lite:general_tac} the general tactical design decisions that arise in a ``standard'' SDD system, wherein orders placed by a cutoff time are delivered by vehicles that return to the depot by the deadline. Then, we will present in Section \ref{lite:newVar} some new variants arising in SDD associated with various features, such as fairness, order cutoff time, and more.} The papers are summarized in Table \ref{table:literature_district} including the first two papers on general tactical decisions, followed by a paper on service network design with hub capacity, and the last one on issues related to SDD service rates. 
\begin{table}[!h]
\caption{Summary of papers on districting. The abbreviations: single-vehicle (SV), multiple trips per vehicle (MT), multiple depots (MD), vehicle capacity (VC), time windows (TW), delivery deadlines (DD), and reinforcement learning (RL).}
\label{table:literature_district}
\scalebox{0.62}{\color{black}
\begin{tabular}{lp{0.18\textwidth}p{0.2\textwidth}p{0.2\textwidth}p{0.12\textwidth}p{0.2\textwidth}lp{0.12\textwidth}}\hline
\textbf{Literature}     & \textbf{Problem}    & \textbf{Objective}                                                                          & \textbf{Uncertainty}                                                                       & \textbf{Uncertainty Nature} & \textbf{Decision}              & \textbf{Constraint} & \textbf{Solution Approach}            \\ \hline
\cite{stroh2022tactical}    & tactical design of SDD systems               & min(delivery time)    &  requests(arrival time, location)    & dynamic \& stochastic  & fleet sizing, length of service window, and consolidation with overnight orders & DD; VC; MT; SV & distribution modeling approach; continuous approximation                   \\
\cite{banerjee2022fleet} & tactical design of SDD systems 
& min(\#vehicle to serve the region)                            & requests(arrival time)    & dynamic \& stochastic                  & fleet sizing, service region partitioning                                                                  & DD; MT     & dispatching policies with Voronoi partitioning; continuous approximation \\
\cite{yildiz2019service} & service and capacity planning in on-demand meal delivery platforms & max(profit of an online meal delivery platform)                            & requests(arrival time), mean delivery time, courier(arrival time), delivery offer acceptance-probability for couriers    & dynamic \& stochastic                  & delivery capacity and composition, delivery service charge, courier compensation, service area coverage                                                                  & DD     & modeling framework with theoretical analysis and numerical examples \\
\cite{chen2022same}       & VRP                                          & max(overall service level and minimum regional service level) & requests(arrival time, location)                                                                                                             & dynamic  \& stochastic & service region sizing; accept/reject; route                                                                                      & DD; MT & RL                                                  \\ 
\cite{banerjee2023has}       & tactical design of SDD systems                                          & max(\#number of requests served) & requests(arrival time, location)                                                                                                             & dynamic \& stochastic & order cutoff times, service area coverage                                                                                      & DD; MT & nonlinear programming with continuous approximation                                                  \\ 
 \cite{wu2022service}    & service network design for SDD system        & min(cost of a vehicle schedule that maximizes \#commodities for which a feasible path can be found) & no & static \& deterministic                   & route; schedule                                                                                            & TW; DD; VC; MD & three heuristics                                                 
                                
\\\hline
\end{tabular}}
\end{table}

\subsection{Literature review on general tactical design decisions}
\label{lite:general_tac}

\cite{stroh2022tactical} propose a vehicle-dispatching model for answering some tactical design questions in SDD distribution. Their model uses continuous approximation techniques to capture the average behavior of an SDD system. Specifically, they consider deliveries from a single {distribution center}, and major decisions are concerned with selecting the optimal
fleet size, determining an order cutoff time, and combining SDD with overnight order delivery operations. {The objective is to minimize the total delivery time spent by all vehicles in serving the orders received within a cutoff time.} {The computational study  using realistic data shows the model is good at predicting the system behavior.} 

\cite{banerjee2022fleet} also use continuous approximations to capture average-case operational behavior. However, {rather than assuming every vehicle may serve orders from the whole region as in \cite{stroh2022tactical},} this paper focuses on deciding the minimum number of vehicles to serve a service region when it can be partitioned into zones dedicated to individual vehicles. The proposed method  first maximizes the area associated with a single-vehicle delivery zone independently and then determines optimal areas as a function of the distance from the depot. As for the first step, the authors propose a Voronoi approach to partition the service region into single-vehicle zones. The Voronoi approach creates zones by assigning to each generator location in a Euclidean plane 
the set of delivery locations closer to it than to any other generator location. Afterwards, the fleet size is estimated based on the surface of the service area. 
The simulations suggest that combining the proposed simple operational dispatching policy with the Voronoi partitioning produces solutions that meet order cutoff time in most cases.

\subsection{Literature review on new variants}
\label{lite:newVar}
Based on traditional logistic districting problems focusing on fleet sizing and routing, several new variants associated with, e.g., planning crowdsourced delivery capacity, customer fairness, fleet assignment with heterogeneous fleets, assignment of couriers on-demand, and hub capacity, appeared in the context of SDD. In the following, we elaborate on these variants and present research works dedicated to each of them. 

\textit{Planning crowdsourced delivery capacity:} \cite{yildiz2019service} investigate service and capacity planning problems in crowdsourced meal delivery. Unlike \cite{stroh2022tactical} and \cite{banerjee2022fleet} who focus on fleet size optimization or maximization of vehicle utilization for an SDD system, their model aims to maximize the profit of an online meal delivery platform while ensuring a target service quality level. They analyze factors like service area coverage, courier compensation, probabilities of courier order acceptance, and composition of delivery capacity (crowd-sourced couriers and company-employed drivers). The analysis with numerical examples offers several interesting insights, including the finding that the hybrid delivery capacity (crowdsourced couriers and company-employed drivers) can improve both reliability and profit.

\textit{Customer fairness:} Concerns about the fairness of the delivery services are rising as fast delivery services are spreading out. In 2016, Amazon was accused of racial disparities in its Prime SDD zones in several U.S. cities (\cite{Amazon_fair}). To address the unfairness in SDD service availability, \cite{chen2022same} first partition the service area and then
{define} fairness as the minimum service rate over the resulting regions.  They model the problem as a multi-objective MDP with a weighted sum of two objectives: maximizing the overall service level (utility) and maximizing the minimum service rate across all regions (fairness), controlled by a factor ranging from 0 to 1. When the factor is valued at 0, no fairness is considered; when it is at 1, no utility is considered. They develop a deep Q-learning approach and show its effectiveness in moderating unfairness for customers from different geographic areas. 

The work of \cite{banerjee2023has} also contribute to the fairness literature. Unlike the papers \citep{stroh2022tactical, banerjee2022fleet, chen2022same} mentioned above, that assume that orders arrive from a predefined service region before a cutoff time, \cite{banerjee2023has} consider the service region as a decision variable and allow it to vary throughout the day. Their goal is to select a service region and deadline that maximize the expected number of orders served. 
They offer various order cutoff times to different areas. The authors find that variable-area systems yield higher revenue. However, fixed-area systems offer a higher level of fairness to customers within the systems.

\textit{Fleet assignment with heterogeneous fleets:} Another possibility to be effective in SDD services is to partition the service region according to the geographical location of requests and decide the delivery mode for each subregion according to the viability of each mode, as we have seen in Section \ref{sec:deliveryMean}. 
For example, \cite{ulmer2018same}, whose work was introduced in Section \ref{sec:deliveryMean}, {propose to divide the service area into two zones, one of which is served by drones, and the other by vehicles. }

\textit{Assignment of couriers on-demand:} In a rapid delivery system, the assignment of couriers to demand points is a crucial decision \citep{agatz2012optimization, savelsbergh201650th, archetti2021recent}. 
The work of \cite{auad2021courier}, introduced in Section \ref{meal_delivery}, {proposes a dynamic couriers assignment policy} 
for SDD. 
Initially, small operating regions are allocated to couriers,  but their boundaries can be enlarged in an on-demand fashion.
That way, flexibility in assignment decisions is introduced, with the aim of increasing couriers' working efficiency. Couriers' satisfaction is improved, as they are no longer expected to wander over the full-service region.
The  computational experiments show that the proposed strategy increases courier satisfaction with a minor impact on customer satisfaction.

\textit{Hub capacity:} 
{Hubs are facilities that switch, sort, connect, and consolidate/disaggregate traffic between multiple origins and destinations (\cite{alumur2021perspectives}).} Limited hub capacities restrict the number of vehicles that can simultaneously be loaded and unloaded in a hub. In urban SDD, with an aggressive amount of orders, hub capacity constraints can lead to waiting at hubs,  which, in turn, renders specific shipment paths  
time infeasible. \cite{wu2022service} investigate the problem with hub capacity constraints. {The goal is to determine a routing path for each commodity so that the loading and unloading times together with capacities at hubs are respected.} 
{They consider a lexicographic objective function that first maximizes the number of commodities served and then minimizes the total travel distance.} They  model the problem on a time-expanded network, as commonly done when time is critical (\cite{belieres2021time}). 
Each time-expanded node stands for a stakeholder, i.e., supplier, warehouse, or customer, at a given period of time. 
The problem is formulated as a MILP and three heuristics are proposed  
to solve real-world instances efficiently. 
This is one of the very few problems that tackles tactical decision problems in SDD, different than districting.

\section{Model Comparison}
\label{model-comparison}
This section analyzes structural differences of optimal solutions obtained when solving formulations \eqref{eq:base1}-\eqref{eq:base4}. The analysis of solutions and illustrative examples given in the following serve several goals:
\begin{itemize}
    \item Concept clarification: they help elucidate and clarify the underlying concepts associated with the problems studied in the literature.
    \item Parameter sensitivity: they help validate if specific parameter values can exert a notable influence on the outcomes.
    \item Impact of objectives and constraints: they help illustrate how solutions change when different constraints or objectives are imposed. 
\end{itemize}  
We compare the obtained solutions in terms of customer fairness, quality of service, and possible environmental impacts. For each problem model, we list metrics related to service quality, fairness, and environmental impact. Specifically, travel distance measures environmental impact, service rate and average waiting time measure service level, and max waiting time and waiting time variability (the difference between max waiting time and min waiting time) measure fairness. Besides, except for the service rate, all the aforementioned metrics are calculated based on the requests served.

\begin{figure}[]
    \centering
    \includegraphics[width=\textwidth]{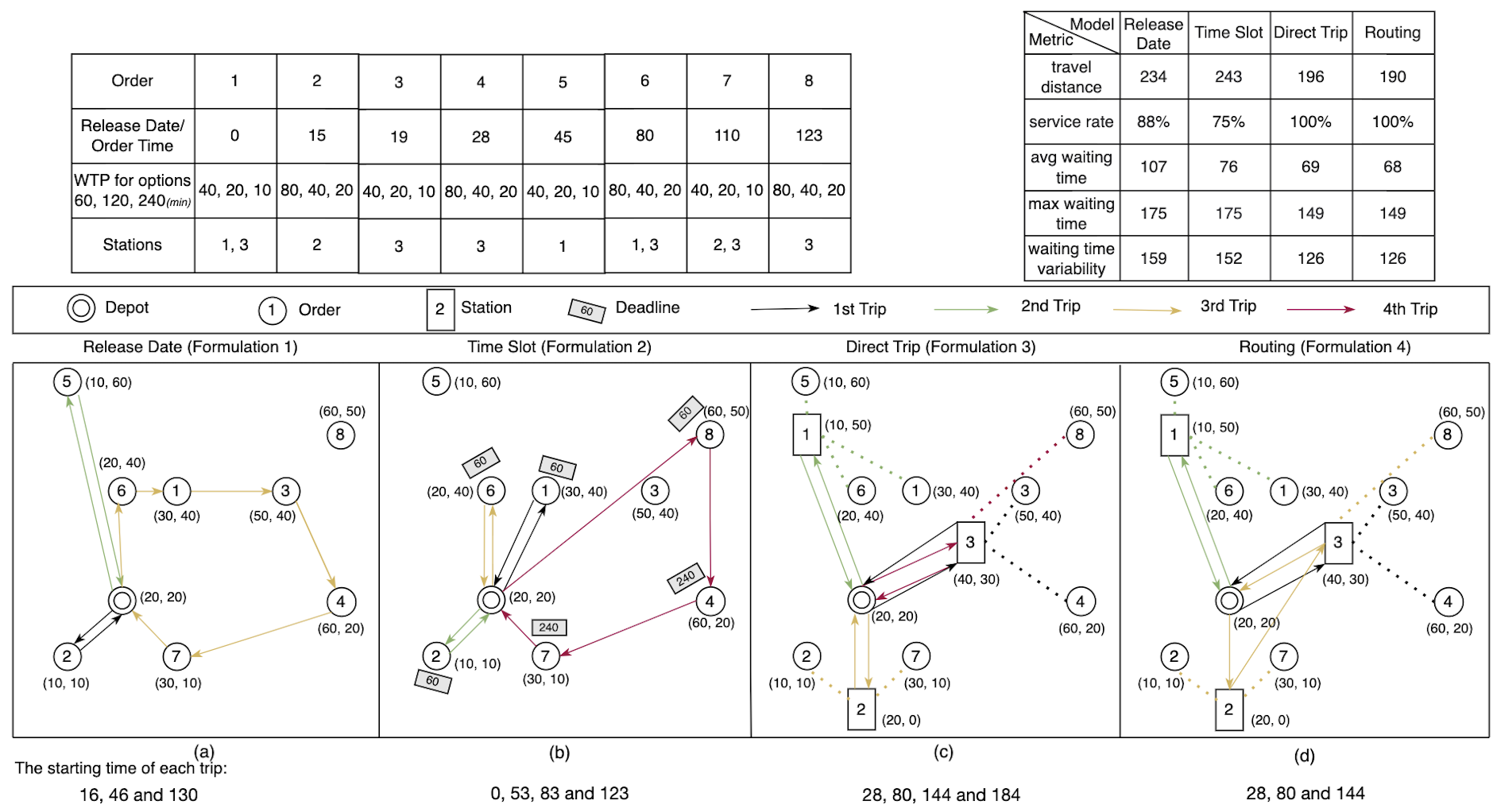}
    \caption{Small instance and corresponding solutions.}
    \label{fig:comparison}
\end{figure}

Data about the first instance and its optimal solutions associated with each of the four formulations are presented in Figure \ref{fig:comparison}.
The instance has eight orders to be delivered within the time horizon [0, 250]. Each order is associated with a release date valued from 0 to 123, which remains the same over all formulations. 
As for the time slot option problem  (Formulation \eqref{eq:base2}), there are three SDD options, i.e., $|\Delta|=3$ {with deadlines $D_\delta \in \{ 60, 120, 240\}$ minutes.  WTP for each customer and each SDD option is provided in the third row of the table shown in Figure \ref{fig:comparison}.}  For example, the customer requesting order 1 is willing to pay at most 40 to get the order delivered within 60 minutes. Concerning the problem with delivery options (Formulation \eqref{eq:base3}), there are three stations with a storing capacity of 3 orders. We determine the set of stations accepted by each customer according to the travel time radius (30 time units) away from the customer's home. The last line of the top left table lists the potential stations for delivering orders. For instance, order 1 can be delivered to station 1 or 3 as long as there is available capacity. 
{The coordinates of points are denoted in brackets in the squares and we use Euclidean distances.} {Circles represent customer locations, whereas rectangles correspond to pickup stations. The central depot is shown as a double circle.} Finally, we fix the value of $M$ to 250 in formulations \eqref{eq:base3} and \eqref{eq:base4}.

\begin{figure}[]
    \centering    \includegraphics[width=0.63\textwidth]{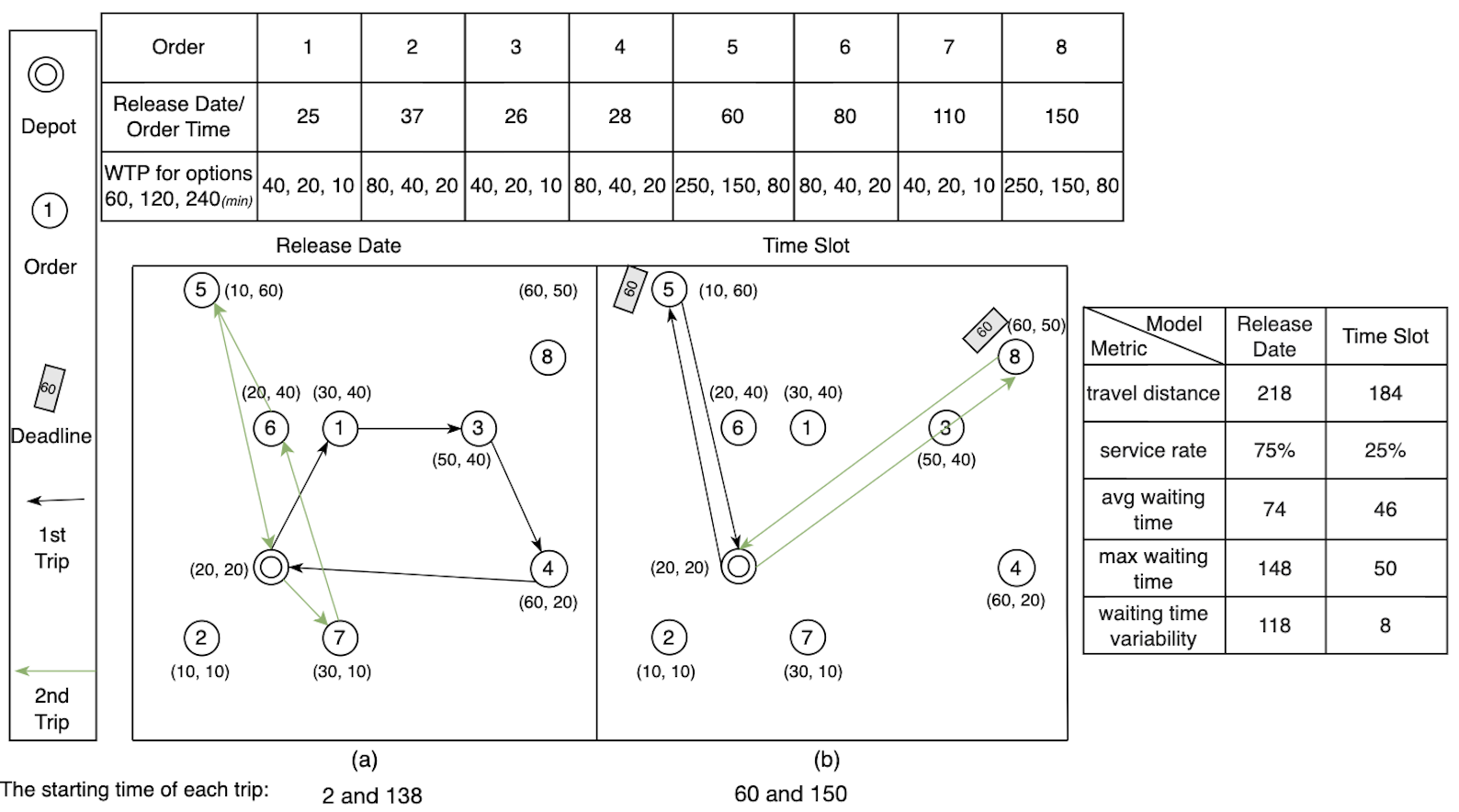}
    \caption{New instance and solutions obtained after changing the WTP and release dates of orders 5 and 8.} 
    \label{fig:extremeWTP}
\end{figure}

The bottom part of Figure \ref{fig:comparison} represents, 
from left to right ((a) to (d)), the optimal solutions obtained by solving formulations \eqref{eq:base1}--\eqref{eq:base4}, respectively. {Figure \ref{fig:comparison}(a)}  shows the solution of the problem with release dates (Formulation \eqref{eq:base1}). We observe that seven requests can be served by three trips before the deadline. 
Comparing it with {Figure \ref{fig:comparison}(b)}, the solution of the time slot problem (Formulation \eqref{eq:base2}), 
we notice that order 8 instead of orders 3 and 5 is served, which is 
due to the fact that order 8 has a higher WTP which leads to a higher profit.  {This example illustrates that due to the maximization of the expected revenue, not only fewer parcels may be served within the same deadline, but also more trips may be needed.}
Indeed, one more trip is dispatched for the solution of Formulation \eqref{eq:base2} compared to the solution of Formulation \eqref{eq:base1}. 
In the time slot problem, 
if there are a few late customer orders with very high WTP located very far from the depot, it might be that the solution serves those few customers and leaves many other customers unserved. This is illustrated in Figure \ref{fig:extremeWTP}(b) where the WTP of orders 5 and 8 has been increased to (250, 150, 80) {and their release dates have been increased from 45 and 123 to 60 and 150 minutes, respectively}. We see that in this case customers 5 and 8 only are served. The solution is compared with the one of Formulation \eqref{eq:base1} (with release dates, depicted in Figure \ref{fig:extremeWTP}(a)), where six requests are served instead. 
This highlights the fact that the time slot model might disappoint certain customers. This is detrimental to retaining customers' loyalty in a highly competitive environment. {Furthermore, the pursuit of maximum profit may obstruct consolidation opportunities and lead to a higher frequency of trips, causing negative externalities on the environment,} which can be seen in the upper-right table of Figure \ref{fig:comparison}: while the time slot model serves the fewest customers, it incurs the greatest travel distance. A possible improvement on the model can be achieved by introducing two hierarchical objectives: maximizing the number of customers served and maximizing the revenue.

\begin{figure}[]
    \centering    \includegraphics[width=\textwidth]{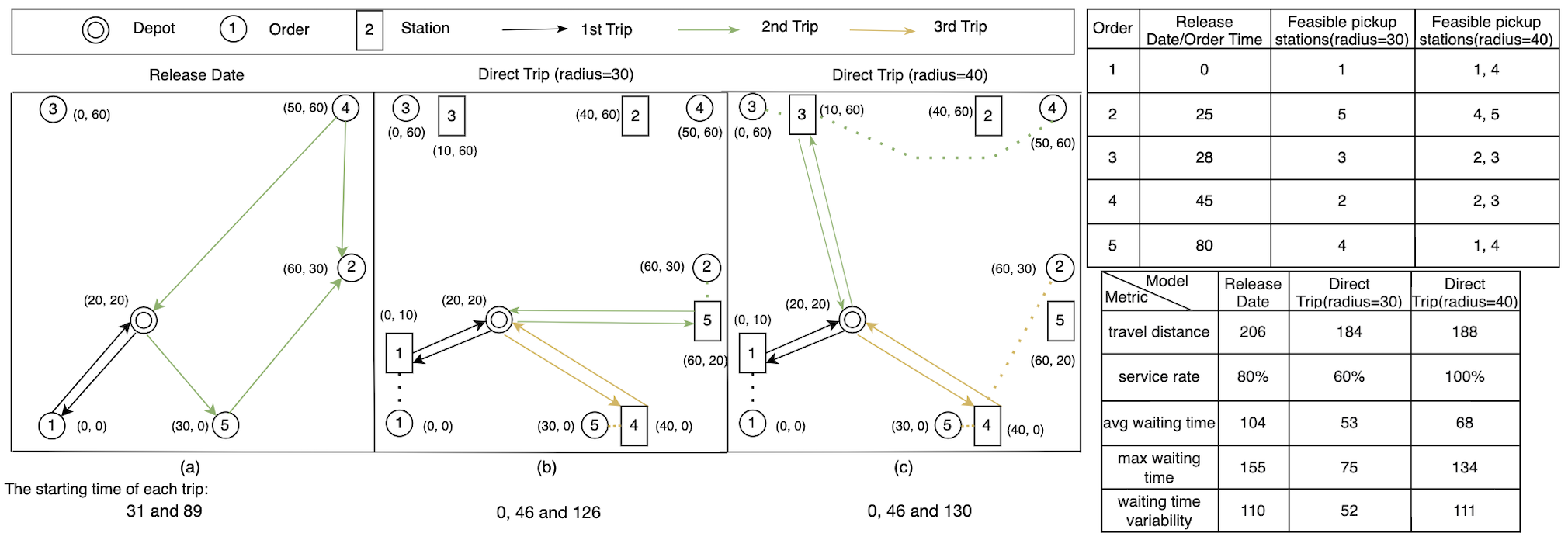}
    \caption{New instance and solutions obtained when changing the radius for determining feasible pickup stations.}
    \label{fig:comparison3}
\end{figure}

Now let us turn our attention to the problems with delivery options. Formulation \eqref{eq:base3} models the case where only direct trips are allowed to replenish pickup stations.
In the solution shown in Figure {\ref{fig:comparison}(c),} four trips serve all requests with a service time of 553. On the other hand, when enabling visiting more than one pickup station in a trip (Formulation \eqref{eq:base4}), only three trips are needed to serve all requests, 
and the service time is reduced to 547. The upper-right table of Figure \ref{fig:comparison} also shows both the direct trip (Formulation \eqref{eq:base3}) and routing models (Formulation \eqref{eq:base4}) get shorter travel distances and reduced average and maximum waiting time when compared to the previous two models, which shows their better performance with respect to environmental impact, fairness, and service rate. Furthermore, the routing model achieves a slightly shorter travel distance and less average waiting time when compared to the direct trip model.  
Nevertheless, direct trip models can still be interesting, if eco-friendly delivery alternatives are used, such as droids or drones (\cite{ulmer2019same}). Indeed, it is often assumed that these two means of transportation can handle one package at a time. 
An extreme situation arises 
when there is only a single feasible pickup station per each customer. 
Figure \ref{fig:comparison3} illustrates such a case: {there are five orders with release dates ranging from 0 to 80 to be delivered before the deadline $T_E=250$.  In the table of Figure \ref{fig:comparison3}, the column named \textit{Stations 1} lists the feasible pickup stations for {formulations \eqref{eq:base3} and \eqref{eq:base4}} with radius equal to 30, and \textit{Stations 2} lists the ones with a bigger coverage, i.e., radius equal to 40. The capacity of each station is 2.} The solution shown in  Figure \ref{fig:comparison3}(b) is obtained by setting the radius equal to 30, in which case, each pickup station can serve one customer only. 
We notice that only three customers can be served before the deadline, and the service time is 586. When increasing the radius from 30 to 40, we observe that all customers can be served (see Figure \ref{fig:comparison3}(c)), and the objective value is reduced to 339.  In the lower-right table of Figure \ref{fig:comparison3}, we see that, in comparison to the previous model, the direct trip model with a bigger coverage improves service rate from 60\% to 100\% with a marginal increase in travel distance and average waiting time but with a bigger increment of max waiting time and waiting time variability.
By comparing the solution obtained using pickup stations with the one obtained using the release dates model (see Figure \ref{fig:comparison3}(a)), we observe that pickup stations may be advantageous in terms of the number of requests served and the total service time. However, this is not always the case and it highly depends on the location of pickup stations and the radius customers are willing to move in order to get their package. {Thus, these two factors have to be carefully considered during strategic and tactical design phases, before deploying pickup-station-based models}. 

To conclude, we observe that changing constraints and objectives can lead to different routes and order assignments for the same instance. To select a fair and robust model, a sensitivity analysis including extreme case instances may help in better understanding the pros and cons of the proposed models. A deliberate selection of objectives, constraints, and broad tests to fulfill both the customer needs and company expectations,  {while keeping sustainability goals in mind,} is imperative.



\section{Challenges and future research directions}
\label{challenges}

Despite the continuous developments of new delivery technologies and advances in analytical methods, providing satisfactory services to customers while making positive profits is still challenging in SDD services. Besides, environmental impacts associated with SDD operations are nonnegligible.
We summarize in the following the most challenging topics that are partially covered in the SDD problems reviewed in this survey, as well as those that still have to be addressed by future research.

\textit{Increased uncertainties coupled with tight delivery times:} There are several sources of uncertainty in SDD, such as unknown customer requests arriving during the day, uncertain customer locations and demand, traffic jams during peak hours, or unexpected road works that may require rerouting. 
Some of these sources of uncertainty affect other distribution problems as well. However, in SDD they have a stronger impact due to the very tight delivery times. In routing problems with release dates and deadlines (see Section \ref{RD_deadlines}), time uncertainty is one of the major obstacles to consolidation. 
For home-attended deliveries (see Section  \ref{delivery}),
in order not to miss the delivery deadline, 
 logistic companies 
 often propose to customers a set of delivery time slots or different delivery options. However, this renders the corresponding optimization problems even more challenging to solve. \cite{hildebrandt2022opportunities} highlight two challenges in solving such stochastic and dynamic problems - searching for a solution in vast solution space and evaluating each solution considering future information changes and actions.
 In this article, we have outlined recent approaches 
 to cope with uncertainties in order to better exploit consolidation opportunities or to increase the number of orders delivered on time. Exploiting historical data to predict release dates or order times, in combination with anticipative algorithms, is one of the widely used approaches to avoid myopic behaviors.
   However, the research in this area is still limited and we believe it will grow along with the increasing coverage of SDD services in the future. 

\textit{Runtime of online algorithms:} In SDD, algorithms developed for operational decisions must be fast enough to ensure prompt decision-making. In decisions involving direct customer interactions (through online platforms, for example), results must be provided within (milli) seconds, thereby avoiding disappointing customers; in the context of vehicle dispatching, decisions may be made within minutes to enable seamless deliveries. Most literature we reviewed underscores the importance of adhering to strict runtime constraints (e.g., \citet{mitrovic2004waiting, van2019delivery,klapp2020request, li2022reinforcement}). However, with increased geographic scales and more customized options offered, such as multiple delivery options and specific customer preferences, problem settings are becoming more complex, and researchers face challenges in developing new algorithms to keep pace with the rising complexity of problems. 

\textit{Autonomous vehicles and drones:} 
For safety reasons, strict environmental and security requirements are imposed on deliveries using autonomous vehicles and drones, and this prevents their deployment at a larger scale\footnote{https://www.abc.net.au/news/2022-09-30/food-delivery-drone-lands-on-power-lines-qld-browns-plains/101489670}. The risk perceived by consumers challenges the societal acceptance of these two transportation means. {According to a survey\footnote{https://morningconsult.com/2022/07/19/drones-safety-consumer-trust-data/} conducted by Morning Consult, 66\% of the respondents worried about air traffic safety due to drone usage.} Moreover, 71\% of the respondents answered they were concerned about personal or data privacy issues related to drone deliveries.

There have been several studies on drones for SDD as discussed in Section \ref{sec:deliveryMean}. However, \cite{TL_industry} states that transportation and logistics companies are not yet ready to use autonomous vehicles and drones for delivering packages to customers' homes. Moreover, their impacts and challenges are currently unclear as there is no mass deployment yet, and the research on the safety impact is just emerging (\cite{tafidis2022safety}). We believe that safety and security requirements should be integrated into the models to alleviate  these concerns. As emphasized by   \cite{poikonen2021future} and \cite{bakach2022robot}, future models should incorporate better representations of the constraints and the new features related to the use of drones or autonomous vehicles, respectively. Moreover, behavioral economics models based on e.g.,  willingness to accept, may be used to address societal concerns.

\textit{Sustainability:} Faster may mean better for customers but not for the planet.  \cite{munoz2021environmental} show that fast deliveries increased total $CO_2$ emissions by up to 15\% as they reduced the possibilities of cargo consolidation. However, the demand for faster shipping is growing, and the pressure on companies is high to avoid losing their hard-earned sales. 
Most scientific studies ignore this issue and focus solely on revenue-related aspects. To have a sustainable future, new delivery strategies need to be envisaged that incorporate environmental aspects. The environmental aspect should also be considered in model objectives to undertake environmental responsibility. To this aim, the first step is measuring the environmental impact of vehicles' operations in terms of emissions (\cite{lin2018demand}). Stakeholders should invest more in improving the efficiency of the delivery network and utilizing greener transportation means (like electric vehicles, bicycles, and porters...) and/or setting up pickup stations. Also, getting customers aware of the topic contributes to sustainable deliveries. A possibility is to provide information while customers make choices for stimulating their eco-consciousness (\cite{gleim2013against}). For example, \cite{agatz2021impact} find that the use of green labels that specify environmentally friendly time slots helps in motivating customers to choose greener time slots that reduce routing distance.

\textit{Fairness and road safety:}
Most papers consider profit-maximizing objectives and completely ignore the security of drivers and  road users. According to a recent study\footnote{https://thesoc.org/wp-content/uploads/2022/05/The-Worst-Mile.pdf}, drivers delivering parcels for Amazon Delivery Service Partners\footnote{The workers in Amazon's delivery system operate as direct employees at one of its facilities, send parcels as gig worker via Amazon's app Flex, or drive for one of Amazon's delivery contractors, called Delivery Service Partners.} were injured at nearly one injury per five full-time-equivalent workers in 2021. Gig workers often spend quite a long time on the road, whether rainy or sunny, and they try to complete as many orders as they can to get the incentives offered by companies. Several risk factors are causing safety issues, such as fatigue, rush to fulfill orders, and harsh weather conditions (\cite{cna}). These road safety conditions affect every delivery driver regardless of the kind of transport used. To lower the risk, companies can implement a set of policies and assistance support. An interesting alternative could be to develop algorithms leading drivers to follow a safer path or making them less stressful while delivering orders on time, for example by taking into account the deadline of orders and other conditions like accident-prone areas, weather, riders' riding distance, and experience (to estimate fatigue level, \cite{fu2022long}).
\citet{yildiz2019provably} propose a model that guarantees minimum payment for drivers, without significantly affecting the service level and other typical performance measures. Balanced and fair distribution of the workload among the drivers is another topic that deserves more attention when developing future SDD services.

Focusing on methodology and 
optimization approaches for SDD services, we envisage the following future research directions:

   \textit{Multiple objectives:} 
    In order to select the most proper model for the application under consideration, the impacts of decisions should be considered from different perspectives, such as customers' fairness, drivers' safety, quality of service, number of customers served, or sustainability. Thus, a multi-objective modeling approach would help decision makers to simultaneously take account of different and sometimes conflicting goals.
    {This is particularly important for developing greener and sustainable models and the respective delivery strategies. }
    While some of the recent studies employ tools of multi-objective optimization \citep{pureza2008waiting,yildiz2019provably,chen2022same,wu2022service}, the full potential of this methodology for SDD services is still to be exploited.
 
    \textit{Modeling uncertainty with two-stage or multi-stage stochastic programming models:} 
    As stated above, uncertainty sources are many in SDD services. Many of the papers dealing with stochastic and dynamic problems propose MDP-based approaches, which are fit for online decisions. However, for offline decisions or when combined with strategic and tactical decisions, the approaches using two-stage or multi-stage stochastic programming (\cite{birge2011introduction}) are another alternative. The computational performance of two-stage models is recently studied in \cite{li2022reinforcement} and \cite{ausseil2022supplier}, and the results show the great potential of these approaches. 
    
    \textit{Modeling SDD services with multiple decision makers:} 
    When problems involve multiple decision makers, one could consider them from the game theory perspective. Various SDD problems could be modeled as simultaneous or sequential games, with collaborative or competitive players. 
    For example, in the problem discussed in Section \ref{formu_pricing}, which is to determine prices of delivery time slots so as to maximize the expected revenue, the provider first decides the price of each slot option for each individual, then each customer selects an option based on his or her WTP. This can be translated into a  sequential two-stage (Stackelberg) game:
     the first decision maker (who acts as a leader) decides the prices, and then the multiple followers (customers) decide which time slot to accept (see e.g., \citet{cerulli2023bilevel}).
     
    \textit{Empirical studies comparing pros and cons of different strategies:}
    In Section \ref{model-comparison}, we compared different SDD problems (strategies) in terms of customer fairness, quality of service, and environmental aspects. While this comparison is conducted on small illustrative instances, 
    more detailed empirical studies and simulations using real-world data need to be done. More critical and comparative studies are indeed necessary in a much broader context of SDD services, in order 
    to better understand the challenges and opportunities of different strategies and 
    provide a more solid base for stakeholders in choosing their preferred strategy. Establishing a common and publicly available benchmark set of instances is of enormous importance for such empirical studies. In Table \ref{table:data_url}, we list available test data provided in the papers surveyed. The literature on the first two entries uses data adapted from \cite{solomon1987algorithms}; the last three on meal delivery are crafted based on real-world data. It is evident that there persists a lack of publicly accessible data sources.
    
    

\begin{table}[!h]
\caption{Available test data of papers surveyed.}
\label{table:data_url}
\scalebox{0.68}{\color{black}
\begin{tabular}{ll}\hline
\textbf{Literature}                                                       & \textbf{Data url}                                                                    \\\hline
\cite{voccia2019same}/\cite{cote2023branch}                                               & \url{http://ir.uiowa.edu/tippie\_pubs/65/}                                   \\
\cite{li2022reinforcement}/\cite{archetti2020dynamic}                              & \url{https://drive.google.com/file/d/1mVw2u1vuvq6NooKbfHlcYPQReBJXqimQ/view} \\
\cite{ausseil2022supplier}                                              & \url{github.com/bstabler/TransportationNetworks/tree/master/chicagoregional} \\
\cite{reyes2018meal}/\cite{yildiz2019provably}/\cite{auad2021courier} & \url{https://github.com/grubhub/mdrplib}\\
\cite{mao2019faster} & \url{https://doi.org/10.1287/msom.2022.1112}
\\\hline                                
\end{tabular}}
\end{table}

\section{Conclusion}
\label{conclusion}
To compete with brick-and-mortar retailers, many e-retailers offer SDD services. Fast delivery services are mainly aimed at compensating for the lack of ``instant gratification'' which affects e-shoppers versus standard shoppers. 
These fast delivery services come with new challenges for distribution and logistic companies.
The paper focuses on the optimization problems in SDD services and analyses the corresponding literature. Particularly, we show how researchers address the dynamic arrival times of orders (Section \ref{RD_deadlines}), decisions about time slot allocation for deliveries (Section \ref{time_window}), delivery options selection (Section \ref{delivery}), emerging problems in pickup and delivery systems (Section \ref{pickup_delivery}) and tactical districting decisions (Section \ref{Districting}). Moreover, we highlight the more representative problems emerging in SDD by providing formulations for deterministic versions and discussing the advantages and disadvantages. 
We also provide a comparison of solution structures obtained by applying the different formulations and discuss the corresponding implications.
Finally, we present challenges in dealing with increased uncertainties, guaranteeing reliable deployment of innovative delivery means, environmental concerns, and road safety. We discuss potential research directions associated with these challenges.

The objectives explored in the literature encompass critical aspects of optimizing service quality and financial performance, and the changes in the objectives or constraints have a big impact on the decisions. The objectives of a company may evolve across different stages of its lifecycle: at the startup stage, the key objective may be maximizing survival probability with limited resources, at the growth and maturity stages, the primary objectives may shift to increasing profitability or maximizing some performance measures. Moreover, for large companies, environmental social responsibility should also be more important. Similarly, the constraints also vary across stages, resulting in different optimization problems to solve.

We conclude by noting that the research field associated with SDD services is evolving very fast, as it happens to the corresponding practical counterpart. Thus, developing interesting research results with a long-term impact in the field is a hard task as new results might become obsolete too soon. Still, the field offers a lot of challenges and opportunities to study and develop novel optimization tools for highly complex systems. 

\section*{Acknowledgement}{%
This work was partially funded by CY Initiative of Excellence (grant ``Investissements d'Avenir" ANR-16-IDEX-0008''). This support is greatly appreciated. The authors also like to thank Debajyoti Biswas for his valuable comments on a preliminary version of the paper.
}



\bibliographystyle{elsarticle-harv} 
\bibliography{refe}

\end{document}